\documentclass[a4paper,12pt]{article}
\usepackage{latexsym}
\usepackage{amsfonts}
\usepackage{amsmath}
\def \E {{\rm E }}
\def \R {{\rm I\!R}}
\def \Z {\mathbb{Z}}
\def \P {\cal P}

\def \L {{\rm I\!L}}

\def \h {{ \tilde{h}}}
\def \B {{ \dot{B}}}
\def \f {\tilde{f}}

\def \fc {\bar{f}}
\def \pr {\partial_{{\tilde{f}}}}
\def \pf2 {{\partial_{f}}}
\def \ph1 {{\partial_{h}}}
\def \p {\partial}
\def \A1 {{\partial}{\cal A}_{n}^{(1)}}
\def \A {{\cal A}}
\def \D {{\rm I\!D}}

\def \sf {\sigma_{f}}

\title{Random fields, large deviations and triviality in quantum field theory. Part II}

\author{ Adnan Aboulalaa\footnotemark[1]}
\date{}

\begin{document}
\vskip2cm
\newtheorem{thm}{Theorem}[section]
\newtheorem{defi}[thm]{Definition}
\newtheorem{prop}[thm]{Proposition}
\newtheorem{Rk}[thm]{Remark}
\newtheorem{lm}{Lemma}[section]
\newtheorem{co}[thm]{Corollary}
\maketitle

\footnotetext[1]{ 
Email: adnan.aboulalaa@polytechnique.org
}
 
\begin{abstract}
The approach developed in the first part of this work, partly based on large deviations, led to the non-existence of interacting scalar fields as strong limits of regularized fields in finite volume and dimensions $d\geq 4$. This second part deals with the weak limit problem, for the particular case of $\varphi^{4}$, according to the 3 cases identified in Part I. In two of theses cases, it is shown that the weak limit is trivial in the sense that the limiting field is identically null. Partial results are obtained for the third case. These results are not incompatible with the already known triviality results that led to a Gaussian free field. As a by product of this study, a rigorous formulation of the principle of the least action for quantum scalar fields is established, together with a set of dynamic field equations that provide explicit expressions of the Schwinger functions.

\end{abstract}

\section{Introduction}
In Part I of this work (\cite{AA1}), the possible existence of the $\varphi_{d}^{4}$ or $P(\varphi)_{d}$ Euclidean field in dimension $d\geq 4$ was studied from a constructive point of view, when these fields are obtained as a limit of regularized fields defined through the measure :
\begin{equation}
\label{Intro1}
d\mu_{n} (\varphi) = \frac{e^{- {\cal A}_{n}(\varphi)}}{Z_{{\cal A}_{n}}} d \mu_{0}(\varphi), \; \; \;   Z_{{\cal A}_{n}}= \E e^{- {\cal A}_{n}(\varphi)}.
\end{equation}
where $\mu_{0}$ is the Gaussian free field measure on ${\cal S}^{'}(\R^{d})$ (which is the reference paths space) and ${\cal A}_{n}$ is a well defined action obtained by regularizing the field and introducing renormalization constants. The subscript $n$ refers to a momentum ultraviolet cutoff. In the case of $\varphi_{d}^{4}$, the regularized action is written as:
\begin{equation}
\label{Intro2}
{\cal A}_{n}(\varphi)= \int_{V}  dx [ g_{n}: \varphi^{4}(x):+ m_{n} : \varphi^{2}(x): + a_{n} : (\partial \varphi)(x)^{2}:  ]
\end{equation}
 where $V\subset \R^{d}$ is a finite volume. The usual mass term $(\delta m_{n})^{2}$ is denoted here by $m_{n}$. A suitable transformation of the action to a mean of an array of random variables enabled to make use of some large deviations techniques and to show that $R_{n}=d\mu_{n}/d\mu_{0} \rightarrow 0$ almost everywhere when $d\geq 4$. This implies the impossibility to obtain a $P(\varphi)_{d}$ field as a setwise (strong) limit of the measures $\mu_{n}$ in a finite volume.

In this paper the issue of existence and triviality of the weak limit of the regularized measures is addressed following the approach of Part I. This study leads also to distinguish 3 cases corresponding to which one of the terms $g_{n}$ (case 1), or $m_{n}$ (case 2) or $a_{n}$ (case 3) is dominant in the action ${\cal A}_{n}$. It will be shown that the limiting field is trivial in the cases °2 and °3 and partial results are also obtained for the case °1. Moreover, the limiting field will be vanishing, i.e. the weak limit of the field measures $\mu_{n}$ is $\delta_{0}$.

The triviality issue of QFT was studied within the lattice approximation framework (see \cite{FFS}, \cite{Callaway} and the references therein) and the important results obtained in the early 1980s by Aizenman \cite{Aizenman} and Frohlich \cite{Frohlich} concern the triviality of the $\varphi^{4}$ for single or two components in dimension $d\geq 5 $ with partial results for the case $d=4$. It was shown that the limiting field of the lattice $\varphi^{4}$ is Gaussian under some natural conditions. A widespread (although not unanimous) belief since then is that these results hold for $d=4$. In the same framework, the remaining case in $d=4$ was treated recently by Aizenman and Dumenil-Copin \cite{A-DC}.

{\it A priori}, we should obtain the same results regardless of the regularization procedure used (lattice or momentum cutoff). Despite the apparent difference between the above-mentioned results on triviality and those obtained in this paper and \cite{AA1}, there is no contradiction between the two. 
In fact, in  \cite{Aizenman}, \cite{Frohlich} and \cite{A-DC} the proof of the free Gaussian character of the limiting field, assumes that the two-point function family $S^{\varepsilon}_{2}(x,y)$ is not vanishing as the lattice spacing goes to 0. This assumption allows the control the $A$ factor in the Aizenman-Frohlich inequality by using an infrared bound result (\cite{FSS}, \cite{Sokal}).

The matter of fact is that $S^{\varepsilon}_{2}$ may vanish as $\varepsilon \rightarrow 0 $, in which case the field will be trivial in the sense that $\varphi\equiv 0$, $\mu_{\infty}$-a.e. Let us briefly show how this is possible. At first, let us mention that a series of formulas giving the Schwinger functions (moments of the interacting field under $\mu_{n}$) are exhibited and used in this work. The first of these formulas gives the following expression of the 2d moment:
\begin{equation}
\label{Intro3}
\E_{n} \varphi(f) \varphi(\tilde{h})= (h,f) - \E_{n} \varphi(f) \ph1 \A_{n}, \; \; {\rm with} \; \tilde{h}=C^{-1} h,
\end{equation}
where $\ph1 \A_{n}$ is the (Malliavin) derivative of $\A_{n}$ in the $h$-direction. For the $\varphi^{4}$ field, (\ref{Intro3}) gives:
\begin{align}
\label{Intro4}
\E_{n} \varphi(f) \varphi(\tilde{f})= & (f,f)  - \E_{n} [4 g_{n} \varphi(f) \int_{V} \varphi_{n}^{3}(x) f_{n}(x) dx + 12 g_{n} c_{n} \varphi(f) \varphi_{n}( f_{n})  \nonumber \\
                    & -  2m_{n}\varphi(f) \varphi_{n}( f_{n}) +  2 a_{n} \varphi(f) \varphi_{n}(\Delta f_{n}) ] 
\end{align}
where we use the notation $f_{n}(x)= \delta_{x}^{n}* f (x)$, with $\delta_{x}^{n}$ a $C^{\infty}$ approximation of $\delta_{x}$. Furthermore, we observe that the Griffiths inequalities established in a lattice context are valid for continuous fields with momentum cutoff, and we have $\E_{n} \varphi(f) \varphi_{n}^{3}(x) f_{n}(x) \geq 0$ for functions $f\geq 0$ which implies that:
\begin{align}
\label{Intro5}
\E_{n} \varphi(f) \varphi(\tilde{f}) \leq & (f,f)  - 12 g_{n} c_{n} \E_{n} \varphi(f) \varphi_{n}( f_{n})  \nonumber \\
                    & -  2m_{n} \E_{n}\varphi(f) \varphi_{n}( f_{n}) +  2 a_{n}  \E_{n}\varphi(f) \varphi_{n}(\Delta f_{n}) ] 
\end{align}
 We also observe that $\E_{n}\varphi(f) \varphi_{n}( f_{n}) \sim \E_{n}\varphi(f)^{2}$ and it is then easy to deduce from (\ref{Intro5}) that $\E_{n}\varphi(f)^{2} \rightarrow 0$ when $m_{n}$  dominates the factors $g_{n} c_{n}$ and $a_{n}$.

Besides the dynamic equations of the fields, a second intermediary result obtained in this study concerns the rigorous formulation and proof of the principle of least action for scalar quantum fields. Roughly speaking, it states that the interacting field measure $\mu_{n}$ is supported in the sets that make minimal the classical actions ${\cal A}_{n}(\chi)$ in some sense. However, to make precise the statement of this principle, the action to be minimized will depend on which of the terms $g_{n}$, $m_{n}$ or $a_{n}$ is predominant, a fact that is naturally expected. This gives rise to a secondary action of the type :
\begin{equation}
\label{Intro6}
{\cal A}^{1}_{n}(\chi)= \lambda_{n} \int_{\Lambda_{ni}}:\chi^{4}(x)):dx+ \alpha_{n} \int_{\Lambda_{ni}} :\chi^{2}(x):dx
                            +\beta_{n} \int_{\Lambda_{ni}}:(\partial\chi)^{2}(x):dx
\end{equation}
where $\lambda_{n}, \alpha_{n}$ and $\beta_{n}$ depend on $g_{n}, m_{n}, a_{n}$. The proof relies on a refinement of the large deviations lower bound of \cite{AA1} and some results on Gaussian fields, in particular the small deviations lower bound of Gaussian fields (Talagrand \cite{Talagrand}) and the Gaussian correlation inequality or at least its weak versions.

In section 2, the notations and settings are presented together with some useful results. In particular the correlation inequalities widely used in the lattice approximation framework are stated for continuous fields with ultraviolet cutoff. While the usage of these inequalities simplify the proof to some extent, they do not seem to be as essential as they are in the lattice approach. The dynamic field equations, which appear to be a quite powerful tool are derived in section 3 as a simple consequence of an integration by parts formula. One could establish a link between one of these equations and a field equation derived in Sokal \cite{Sokal} and used in \cite{BFS} for a treatment of $\varphi^{4}$ in dimension 3. In section 4 we state and prove the principle of least action. In section 5 the issue of the existence and triviality of the weak limits of the interacting field measures is addressed according to the above-mentioned cases.

\section{Notations and preliminaries}

We shall use the notations of Part I. The space-time is $\R^{d}, d\geq 2$ and the path space of the Euclidean fields is $\P={\cal S}^{'}(\R^{d})$, i.e. the Schwartz space of tempered distributions on $\R^{d}$, that will be denoted by $\varphi$. It will be the reference probability space with $\mu_{0}$ the Euclidean free field probability measure, that is, the Gaussian measure on $\P$ whose covariance operator is $C=(-\Delta +1)^{-1}$. If $X$ is a Gaussian random variable, then for $p\geq 2$, $:X^{p}:$ will denote its Wick powers. 
\\
\noindent
{\it Tests functions and regularization:} Let $\delta_{x}$ be the Dirac distribution centered in $x\in \R^{d}$ and $f$ an integrable function on $\R^{d}$; we shall use the regularization :
\begin{equation}
\label{Not1}
f_{n}(x)= \delta_{x}^{n}* f (x) = \int f(y)\delta^{n}_{x}(y) dy,
\end{equation}
and 
\begin{equation}
\label{Not1b}
 f_{nn}(x)= \delta^{n}_{x}* f_{n} (x) = \int f_{n}(y)\delta^{n}_{x}(y) dy
\end{equation}
with
\begin{equation}
\label{Not2}
  \delta^{n}_{x}(y)= (2\pi n)^{d/2}e^{-n^{2}(y-x)^{2}/2}    
\end{equation}
And given a distribution $\varphi$, we obtain the smooth functions $\varphi_{n}$ :
\begin{equation}
\label{Not3}
\varphi_{n}(x) := \varphi(\delta^{n}_{x}), \;\; \; \varphi_{n}(f)= \int \varphi_{n}(x) f(x) dx
\end{equation}
We notice that $\delta^{n}_{x} (y)= \delta^{n}_{y} (x)$ and $ \varphi_{n}(f)= \varphi(f_{n})$.
For a smooth function $f\in {\cal S}={\cal S}(\R^{d}) $, the gradient $\nabla f$ of a function will be also noted $\partial f(x)=\sum_{j=1}^{d} \partial_{j} f(x) e_{j}$, $(e_{j}, j=1, ..., d)$, being the canonical basis of $\R^{d}$ and $ (\nabla f)(x)^{2}=(\partial f(x))^{2}= \sum_{j}^{d} (\partial_{j} f(x))^{2}$. The Laplacian $\Delta f$ will be sometimes noted $ \partial^{2} f(x)=\Delta f(x) = \sum_{j}^{d} \partial_{j}^{2} f(x)$. We consider a fixed finite volume $V\subset \R^{d} $ and when there is no ambiguity we shall simplify the notation of the integrals as:
\begin{equation}
\label{Not4}
\int h = \int_{\R^{d}} h(x) dx, \; \; \;  \int_{V} h = \int_{V} h (x) dx 
\end{equation}
\noindent
{\it Action functional and interaction measure}: For an action functional ${\cal A} : {\cal S}'(\R^{d}) \rightarrow \R $, the interaction measure is:
\begin{equation}
\label{Not5}
 d \mu (\varphi) = \frac{e^{- {\cal A}(\varphi)}}{Z_{{\cal A}}} d \mu_{0}, \;  Z_{{\cal A}}= \E e^{- {\cal A}(\varphi)}.
\end{equation}
The Euclidean quantum $\varphi_{d}^{4}$ field model is given by (\ref{Not5}) and the formal action:
\begin{equation}
\label{Not6}
{\cal A}(\varphi)= \int_{\R^{d}}  dx [ g: \varphi^{4}(x): + \; m : \varphi^{2}(x): + \; a: (\partial \varphi)(x)^{2}:]
\end{equation}
For notational convenience, the mass term, i.e. the factor of the $\varphi^{2}(x)$ part, will be denoted by $m$ instead of the usual and natural $m^{2}$ or $\delta m^{2}$ notations. The basic issue is whether such a field (or measure) exists and is non trivial in some sense. To address it, we are led to consider the desired field (\ref{Not5})-(\ref{Not6}) as being approximated by a well-defined regularized field of the type :
\begin{equation}
\label{Not7}
 d\mu_{n} (\varphi) = \frac{e^{- {\cal A}_{n}(\varphi)}}{Z_{{\cal A}_{n}}} d \mu_{0}(\varphi), \; \; \;   Z_{{\cal A}_{n}}= \E e^{- {\cal A}_{n}(\varphi)}.
\end{equation}
The corresponding regularized action for the $\varphi_{d}^{4}$ field is given by:
\begin{equation}
\label{Not8}
{\cal A}_{n}(\varphi)= \int_{\R^{d}}  dx [ g_{n}: \varphi^{4}_{n}(x):+ m_{n} : \varphi^{2}_{n}(x): + a_{n}: (\partial \varphi_{n})(x)^{2}:]
\end{equation}
where $\varphi_{n}$ is a regularized field, leading to well defined Gaussian random variables $\varphi_{n}(x)$ (under $\mu_{0}$) and well defined expressions $: \varphi^{p}_{n}(x):$. This is done by an ultraviolet regularization, such as a lattice or a momentum cutoff. The renormalization constants $g_{n}, m_{n}, a_{n}$ are to be chosen so that the sequence of measures $\mu_{n}$ has a limit. 
In this paper, we consider the momentum cutoff defined by (\ref{Not2})-(\ref{Not3}) and we work first in a finite volume $V$. More general actions can be considered such as polynomial actions ${\cal A}(\varphi) = \int P(\varphi(x), \partial \varphi(x))$ and non symmetric actions like:
\begin{equation}
\label{Not9}
{\cal A}(\varphi)= \int_{\R^{d}}  dx [ g: \varphi^{4}(x):+ m : \varphi^{2}(x): + a: (\partial \varphi)(x)^{2}: + b \varphi(x) ]
\end{equation}
The expectations of a r.v. $Y$ with respect to the measures $\mu_{n}$ will be denoted by $\E_{n}Y$ (or $\left< Y \right>_{n}$):
\[\E_{n} Y = \E_{\mu_{n}} Y = \int Y(\varphi) d \mu_{n}(\varphi) \]
\noindent
We recall the expression of the covariance kernel (cf. \cite{GJ87}):
\[C(x,y)=2^{-1}(2\pi)^{-(d-1})g(|x-y|) \] 
with :
\begin{equation}
\label{Note9b}
 g(t) = \frac{1}{t^{d-2}}|S^{d-2}|\int_{0}^{\infty} \frac{s^{d-3}e^{-s(1+t^{2}s^{-2})^{1/2}}}{(1+t^{2}s^{-2})^{-1/2}} ds
\end{equation}
The behaviour near $0$ is given by $C(u)\sim K |u|^{-(d-2)}$ and a key number of the theory is :
\begin{equation}
\label{Not10}
c_{n} =  \E \varphi^{2}_{n}(x-y) = \E \varphi^{2}_{n}(0) = c_{n} \sim K n^{d-2}
\end{equation}
We also recall an estimate of $c_{n}(x)$ that will be used in this paper:
\begin{align}
\label{Not10b}
c_{n}(x-y) =  \E \varphi_{n}(x)\varphi_{n}(y) &=  \E \varphi_{n}(0)\varphi_{n}(x-y)\sim K\int_{\R^{d}} C(\frac{u}{n}) e^{-(u+(nx-y))^{2}} \\ \nonumber
                                        & \sim  K n^{d-2}\int_{\R^{d}} \frac{e^{-v^{2}}}{|v-n(x-y)|^{d-2}}
\end{align}
See Part I (\cite{AA1}) for estimates of quantities involving the covariances of $\varphi^{p}_{n}(x)$. For $f\in {\cal S}(\R^{d})$ we set 
\begin{equation}
\label{Not11}
\f= C^{-1} f,  \;\;  \fc= C f,  \; \; \;  \sigma_{f}=(f,\f)=(f,C^{-1}f)
\end{equation}

\noindent
For $f^{(1)}, ... f^{(p)} \in {\cal S}$, the Schwinger distributions of the cutoff interaction field are :
\begin{equation}
\label{Not11}
S^{(p)}_{n}(f^{(1)}, ... f^{(p)})= \E_{n} \varphi(f^{(1)}) ... \varphi(f^{(p)})
\end{equation}
The following distribution version of the uniform boundedness principle (see, e.g., Horvath \cite{Horvath}) will be frequently used throughout this paper:
\begin{thm}
\label{Distributions}
For a given $p\geq 1$, suppose that the quantities $S^{(p)}_{n}(f^{(1)}, ... f^{(p)})$ converge some $S^{(p)}_{\infty}(f^{(1)}, ... f^{(p)})$ as $n\rightarrow \infty$ for all $f^{(1)}, ... f^{(p)} \in {\cal S}$ or ${\cal D}$. Then $S^{(p)}_{\infty}$ is a distribution and we also have :
\begin{equation}
\label{Not12}
S^{(p)}_{n}(f^{(1)}_{n}, ... f^{(p)}_{n}) \rightarrow S^{(p)}_{\infty}(f^{(1)}, ... f^{(p)}) \; {\rm as} \;  n\rightarrow \infty
\end{equation}
for all sequences of $f^{(1)}_{n}, ... f^{(p)}_{n}$ that converge to $f^{(1)}, ... f^{(p)}$ in ${\cal S}$ or ${\cal D}$. 
\end{thm}
This fact will be typically used in the following situation : suppose that $\E_{n} \varphi(f)\varphi(h) \rightarrow \E_{\infty} \varphi(f)\varphi(h)$ for all $f, h\in {\cal S}$ then we will also have:
\begin{equation}
\label{Not13}
\E_{n} \varphi_{n}(f)\varphi_{n}(h) \rightarrow \E_{\infty} \varphi(f)\varphi(h), \; {\rm and} \; \; \E_{n} \varphi_{n}(f_{n})\varphi_{n}(h_{n}) \rightarrow \E_{\infty} \varphi(f)\varphi(h)
\end{equation}

\noindent
$\bullet$ {\bf Correlation inequalities for continuous fields with ultraviolet cutoff}
\\
The topic of correlation inequalities has been developed in the context of lattice systems and was also extended to the continuous field when these fields exist and can be approximated by lattice fields (e.g. the case of $P(\varphi)_{2}$, see \cite{GRS75}, \cite{Simon} pp. 278-284, \cite{GJ87}, p. 230). We shall need some of these inequalities for continuous field like $ \varphi^{4}_{d}, P(\varphi)_{d}$ {\it with momentum cutoff}. Let us first recall the kind of inequalities needed.

Let us consider a lattice $\L$, which is a finite set, and a $\varphi = (\varphi_{j}), j\in \L$ a field on $\L$ which. So the $\varphi$ are just functions from $\L$ to $\R$. For a given $A\subset \L$ we use the notation:
\[ \varphi^{A}=\prod_{j\in A} \varphi_{j}  \]
and it will be understood that the sets $A$ may contain multiple copies of a given element $j$. We suppose that we are given a Hamiltonian
\[ H = - \sum_{A \in {\cal P}(\L)}J_{A} \varphi^{A}  \]
where ${\cal P}(\L)$ is the partition set of $\L$. We also suppose that we are given a probability measure $d\nu_{j}(\varphi_{j})$ on $\R$ for each $j$. Then we define the probability measure on the path space of the lattice field $\Sigma:=\R^{\L}$ :
\begin{equation}
\label{lattice1}
 d \nu_{H} (\varphi) = e^{- H(\varphi)} d\nu(\varphi) = e^{- H(\varphi)}  \Pi_{j\in \L}d\nu_{j}(\varphi_{j}) 
\end{equation}
The corresponding expectation will be denoted by $\left<F\right> = \int f(\varphi) d \nu_{H} (\varphi) $. The $\varphi^{4}$ typical example for lattice fields, is given by: 
\begin{equation}
\label{lattice2}
H = - \sum_{i, j} J_{ij}\varphi_{i}\varphi_{j} - \sum_{i} h_{i}\varphi_{i}
\end{equation}
and :
\begin{equation}
\label{lattice3}
d\nu_{j}(\varphi_{j}) = e^{- P_{j}(\varphi_{j})} d \varphi_{j}, \; {\rm with:} \; P_{j}(\varphi_{j})= \lambda_{j} \varphi_{j}^{4} + m_{j} \varphi_{j}^{2}
\end{equation}
Most of the correlation inequalities have been established for ferromagnetic systems, i.e. when $ J_{A} \geq 0$:

$\bullet$ {\it Griffiths inequalities}: If the Hamiltonian $H$ is ferromagnetic ($J_{A} \geq 0$) and the measures $d\nu_{j}(\varphi_{j})$ are symmetric under the map $\varphi_{j} \rightarrow - \varphi_{j}$ then we have:
\begin{align}
\label{Griffiths}
\text{ 1st Griffiths inequality:} \;   & \left<\varphi^{A}\right> \geq 0 \; \forall A \subset \L \nonumber  \\
\text{ 1st Griffiths inequality:} \; &  \left<\varphi^{A} \varphi^{B}\right> \geq  \left<\varphi^{A} \right> \left<\varphi^{B}\right> \; \forall A, B \subset \L
\end{align}
The following inequalities are valid in particular for the $\varphi^{4}$ symmetric systems $h=0$:

$\bullet$ {\it The Lebowitz inequality} :
\begin{equation}
\label{Lebowitz}
\left<\varphi^{i_{1}} \varphi^{i_{2}} \varphi^{i_{3}} \varphi^{i_{4}}\right>  \leq \left<\varphi^{i_{1}} \varphi^{i_{2}}\right> \left<\varphi^{i_{3}} \varphi^{i_{4}}\right> + \left<\varphi^{i_{1}} \varphi^{i_{3}}\right>\left<\varphi^{i_{2}} \varphi^{i_{4}}\right>+ \left<\varphi^{i_{1}} \varphi^{i_{4}}\right>\left<\varphi^{i_{2}} \varphi^{i_{3}}\right>
\end{equation}
This result was generalized to :
\\
\noindent
$\bullet$ The Gaussian inequality (cf. \cite{BFS83} and the references therein):
\begin{equation}
\label{Gaussian}
\left<\varphi^{i_{1}} ... \varphi^{i_{2p}}\right>  \leq \left<\varphi^{i_{1}} ... \varphi^{i_{2p}}\right>^{G} :=  \sum_{pairings \; i_{\alpha},i_{\beta}}\prod_{i_{\alpha},i_{\beta}} \left<\varphi^{i_{\alpha}} \varphi^{i_{\beta}}\right> 
\end{equation}
The notation $\left< ... \right>^{G}$ in right hand side of (\ref{Gaussian}) refers to an expression that has the same structure as the expansion of the expectation of $2p$ Gaussian random variable (but the expectation of two rv's $\left<\varphi^{i_{\alpha}} \varphi^{i_{\beta}}\right> $ is taken with respect to the lattice measure $\nu_{H}$ and not w.r.t a free Gaussian measure), see \cite{BFS83}.
\\
\noindent
$\bullet$ {\it The skeleton inequalities} (\cite{BFS83}) provide upper and lower bounds for the Ursell functions:
\begin{align}
\label{Ursell}
U_{4}:= \left<\varphi^{i_{1}} \varphi^{i_{2}} \varphi^{i_{3}} \varphi^{i_{4}}\right>^{T} & :=  \left<\varphi^{i_{1}} \varphi^{i_{2}} \varphi^{i_{3}} \varphi^{i_{4}}\right> \\ \nonumber
                           & - \left<\varphi^{i_{1}} \varphi^{i_{2}}\right> \left<\varphi^{i_{3}} \varphi^{i_{4}}\right> - \left<\varphi^{i_{1}} \varphi^{i_{3}}\right>\left<\varphi^{i_{2}} \varphi^{i_{4}}\right>
														-  \left<\varphi^{i_{1}} \varphi^{i_{4}}\right>\left<\varphi^{i_{2}} \varphi^{i_{3}}\right>
\end{align}
\begin{equation}
\label{skeleton}
0 \leq -U_{4}  \leq 4! \lambda \sum_{j\in \L} \left<\varphi^{i_{1}} \varphi^{j}\right>\left<\varphi^{i_{2}} \varphi^{j}\right>\left<\varphi^{i_{3}} \varphi^{j}\right>\left<\varphi^{i_{4}} \varphi^{j}\right> 
\end{equation}
The coupling constant $\lambda$ is supposed to be independent of the sites. These inequalities were later extended to other situations (Polynomial interaction \cite{Hara}, higher orders and 2-components systems \cite{BovierFelder}).
\\
\noindent
In our context we have the following:
\begin{prop}
\label{CorrelationIneq}
Let $\mu_{k}, \E_{k} $ be the measure and expectation corresponding to the $\varphi_{d}^{4}$ field in finite volume and with an ultraviolet momentum cutoff $k$, and $(f_{i})$ a sequence of {\it positive functions} on $\R^{d}$. The correlation inequalities are valid for this continuous field with the following expressions:
\\
\noindent
$\bullet$ 1st Griffiths inequality: for $p\geq 2$: 
\begin{equation}
\label{Griffiths1}
 \E_{k} \varphi (h_{1}) ... \varphi(h_{p})  \geq 0 
\end{equation}
$\bullet$ 1d Griffiths inequality: for $p, q \geq 1, p\leq q-1$:
\begin{equation}
\label{Griffiths2}
                   \E_{k} \varphi (h_{1}) ... \varphi(h_{q}) \geq  \E_{k} \varphi (h_{1}) ... \varphi(h_{p}) \E_{k} \varphi (h_{p+1}) ... \varphi(h_{q}).
\end{equation}
$\bullet$ Gaussian inequalities : for $p\geq 1$: 
\begin{equation}
\label{GaussianIneq}
\E_{k} \varphi (h_{i_{1}}) ... \varphi (h_{i_{2p}})  \leq \E_{k}^{G} \varphi(h_{i_{1}}) ... \varphi(h_{i_{2p}}) 
                      :=  \sum_{pairings \; i_{\alpha},i_{\beta}}\prod_{i_{\alpha},i_{\beta}} \E_{k} \varphi (h_{i_{\alpha}}) \varphi(h_{i_{\beta}})
\end{equation}
\end{prop}
{\it Sketch of the proof.} The measure and action under consideration are given by (\ref{Not7}) and (\ref{Not8}) and the cutoff $k$ is now fixed  once for all. The idea is to approximate the measure $\mu_{k}$ by a lattice measure and to notice that it is actually a ferromagnetic or a $\tilde{\varphi}^{4}$ lattice measure. The correlation inequalities are valid in this case and one has to check that they remain valid when we take linear combinations of the lattice variables $\tilde{\varphi}_{j}$. The passage to the continuous limit when the lattice spacing goes to $0$ will give the correlation inequalities of the proposition. This limit is not problematic here because the cutoff $k$ is fixed.

So let us consider the lattice $\L_{\varepsilon}=(\varepsilon \Z)^{d}$, with $\varepsilon> 0$ and $\L_{\varepsilon, V}=\L_{\varepsilon}\bigcap V$. As $d \mu_{k} = \exp(-\A_{k})/Z_{k} d\mu_{0}$, we have to approximate both $\mu_{0}$ and $\A_{k}(\varphi)$. The action (\ref{Not8}) may be approximated first by a Riemann sum of the form:
\begin{equation}
\label{Correl1}
{\cal A}_{k, \varepsilon}(\varphi)= \frac{1}{N_{\varepsilon}} \sum_{j=1}^{N_{\varepsilon}} [g_{k}: \varphi^{4}_{k}(x_{j}):+ m_{k} : \varphi^{2}_{k}(x_{j}): + a_{k}: (\partial \varphi_{k})(x_{j})^{2}:]
\end{equation}
We have to make two modifications to this sum : first the terms $\varphi_{k}(x_{j})$ should be replaced by terms of $\varphi (x_{j})$ (without the $k$ because the measure $\mu_{0}$ depends on $\varphi$ - that will be supposed to be continuous in the lattice approximation, instead of being in ${\cal S}^{'}$, see below). We write :
\begin{equation}
\label{Correl2}
\varphi_{k}(x_{j})=  \int \varphi(y)\delta^{k}_{x_{j}}(y) dy \sim \frac{1}{N_{\varepsilon}} \sum_{l=1}^{N_{\varepsilon}}  \varphi(x_{l})\delta^{k}_{x_{j}}(x_{l})
\end{equation}
In (\ref{Correl2}) and in the following the symbol $\sim$ means that we have an equality when the limit $\varepsilon \rightarrow 0$ is taken in the expression considered.
The second modification of the discrete action (\ref{Correl2}) consists in replacing the gradient term by a finite difference. We have $(\partial \varphi_{k})(x_{j})^{2}=\sum_{i=1}^{d} (\partial_{i} \varphi_{k})(x_{j})^{2}$ and for each $i=1, ..., d$ we take:
\begin{equation}
\label{Correl3}
\partial_{i} \varphi_{k}(x_{j}) \sim \frac{\varphi_{k}(x_{j}+\varepsilon e_{i}) - \varphi_{k}(x_{j}) }{\varepsilon}
\end{equation}
where $e_{i}, i=1, ..., d $ is the  canonical basis of $\R^{d}$. With the modifications (\ref{Correl2}) and (\ref{Correl3}) we get a discrete action of the form:
\begin{equation}
\label{Correl4}
{\cal A}_{k, \varepsilon}(\varphi)= \frac{1}{N_{\varepsilon}}  \sum_{j=1}^{N_{\varepsilon}} (g_{k} \varphi^{4}(x_{j})+ m_{k}^{'} \varphi^{2}(x_{j})) 
                               - J_{k} \sum_{j, j'=1, j\neq j'}^{N_{\varepsilon}} \varphi^{2}(x_{j}) \varphi^{2}(x_{j'})) + C(k)
\end{equation}
where $J_{k}=2a_{k}\varepsilon^{-2}$ and $m_{k}^{'}, C(k)$ are constants depending on $k$. By this approximation, we obtain a first partially discrete measure :
\begin{equation}
\label{Correl5}
 d\mu_{k,\varepsilon}^{(1)} (\varphi) = \frac{e^{- {\cal A}_{k, \varepsilon}(\varphi)}}{Z_{k,\varepsilon}} d \mu_{0}, \;  Z_{k,\varepsilon}^{(1)}= \E e^{- {\cal A}_{k, \varepsilon}(\varphi)}.
\end{equation}
Now, let us consider the lattice $\L_{V,\varepsilon}=(\varepsilon \Z)^{d}\bigcap V$ and for $j=(j_{1}, ..., j_{d})\in \Z^{d}$ we set $x_{j}=\varepsilon j$. The lattice paths are the functions :
\begin{equation}
\label{Correl6}
 \tilde{\varphi}: \L_{V,\varepsilon} \rightarrow \R :  \; x_{j}\mapsto  \tilde{\varphi}_{j}\equiv  \tilde{\varphi}(x_{j})
\end{equation}
The space of such discrete paths or functions is denoted $l^{2}(\L_{V,\varepsilon})$ (we add a square-integrability condition : $\sum_{j}\tilde{\varphi}(x_{j})^{2} <0 $). Now, to the action ${\cal A}_{k, \varepsilon}(\varphi)$ on the space of continuous function $C(V)\subset {\cal S}^{'}$ we associate the action:
\begin{equation}
\label{Correl7}
\tilde{\cal A}_{k, \varepsilon}(\tilde{\varphi})= \frac{1}{N_{\varepsilon}}  \sum_{j=1}^{N_{\varepsilon}} (g_{k} \tilde{\varphi}^{4}(x_{j})+ m_{k}^{'} \tilde{\varphi}^{2}(x_{j})) 
                               - J_{k} \sum_{j, j'=1, j\neq j'}^{N_{\varepsilon}} \varphi(x_{j}) \varphi^{2}(x_{j'})) + C(k)
\end{equation}
defined on $l^{2}(\L_{V,\varepsilon})$. In order to a get a full lattice measure approaching $\mu_{k}$, it remains to approximate the free field Gaussian measure $\mu_{0}$. Let $\mu_{0, \varepsilon}$ be such a measure defined on $l^{2}(\L_{V,\varepsilon})$, see \cite{GJ87}, $\S$ 9.5-9.6 and \cite{Simon}. Our lattice measure is then:
\begin{equation}
\label{Correl8}
d\tilde{\mu}_{k,\varepsilon}(\tilde{\varphi}) = e^{-\tilde{\cal A}_{k, \varepsilon}(\tilde{\varphi})} d\mu_{0} (\tilde{\varphi})
       = e^{-\tilde{\cal A}_{k, \varepsilon}^{(1)}(\tilde{\varphi})} \Pi_{j} d\tilde{\varphi}_{j}
\end{equation}
is a ferromagnetic and $\varphi^{4}$ lattice measure for which the previous correlation inequalities are valid. Also, this measure induces a measure on ${\cal S}^{'}(\R^{d})$ and we have for instance the Gaussian inequality in the lattice form
\begin{equation}
\label{Correl9}
\E_{\tilde{\mu}_{k,\varepsilon}} \tilde{\varphi}_{1} ... \tilde{\varphi}_{2p}  \leq \E_{\tilde{\mu}_{k,\varepsilon}}^{G} \tilde{\varphi}_{1} ... \tilde{\varphi}_{2p}
\end{equation}
The fact that enables the passage from this form to the continuous measure form where we deal with continuous fields smeared with positive functions, is that (\ref{Correl9}) is stable under linear transformation with positive coefficients : if $\tilde{\varphi}_{j}^{'}=\sum_{l} \alpha_{l}\tilde{\varphi}_{l}$ with $\alpha_{l}\geq 0$, we still have:
\begin{equation}
\label{Correl10}
\E_{\tilde{\mu}_{k,\varepsilon}} \tilde{\varphi}_{1}^{'} ... \tilde{\varphi}_{2p}^{'}  \leq \E_{\tilde{\mu}_{k,\varepsilon}}^{G} \tilde{\varphi}_{1}^{'} ... \tilde{\varphi}_{2p}^{'}
\end{equation}
hence, if $f_{q}$ are positive functions on $\R^{d}$, we consider the Riemann sums:
\begin{equation}
\label{Correl11}
\tilde{\varphi}(f_{\varepsilon}):= \frac{1}{N_{\varepsilon}} \sum_{l=1}^{N_{\varepsilon}}  \tilde{\varphi}_{l}f(x_{l}) 
                           \equiv \frac{1}{N_{\varepsilon}} \sum_{l=1}^{N_{\varepsilon}} \varphi(x_{l})f(x_{l})   \equiv  \varphi(f_{\varepsilon})
\end{equation}
(To a continuous function $\varphi$, we associate the lattice function $\tilde{\varphi}$ defined by $\tilde{\varphi}(x_{j})= \varphi(x_{j})$). From (\ref{Correl10}) we deduce that for every positive functions $f_{1}, ..., f_{2p}$ we have :
\begin{equation}
\label{Correl12}
\E_{\tilde{\mu}_{k,\varepsilon}} \tilde{\varphi}(f_{1, \varepsilon}) ... \tilde{\varphi}(f_{2p, \varepsilon})  \leq 
                \E_{\tilde{\mu}_{k,\varepsilon}}^{G} \tilde{\varphi}_{1}(f_{1, \varepsilon}) ... \tilde{\varphi}_{2p}(f_{2p, \varepsilon})
\end{equation}
and when the limit as $\varepsilon \rightarrow 0$ is taken we will have:
\begin{equation}
\label{Correl12b}
\E_{k} \varphi (f_{1}) ... \varphi (f_{2p})  \leq \E_{k}^{G} \varphi (f_{1}) ... \varphi (f_{2p})
\end{equation}
The other correlation inequalities can be proved in the same way. $\Box$

\section{Dynamic equations of the interacting fields}
In this section we shall derive some useful formulas for the moments of random fields defined by an exponential measure of the type (\ref{Not5}). Some of these formulas can be derived directly from the generating functional of the measure and all of them are in fact a consequence of an integration by parts formulas.
\subsection{Integration by parts and dynamic equations}
 
For $h \in {\cal S}(\R^{d})$ we consider the translation on ${\cal S}'(\R^{d})$ : $\tau_{h}: \varphi \rightarrow \varphi+h $ and the induced measure $\mu_{h}$ defined by 
$\mu_{h} (B)= \mu (B) (\tau_{h}^{-1}(B))$ for $B\in {\cal B}$. Then we have :

\begin{prop}
\label{PropQuasiinvariance}
The measures $\mu_{h}$ and $\mu_{0}$ are equivalent and we have :

\begin{equation}
\label{Quasiinvariance1}
 \frac{d \mu_{h}}{d\mu_{0}} = K_{h}(\varphi):= e^{[(\varphi, C^{-1}h) + (h, C^{-1}h)/2]}
\end{equation}
\end{prop}
{\it Proof.} See Glimm-Jaffe \cite{GJ87} p.207.$\Box$
\\
\noindent
For a functional $Y: {\cal S}'(\R^{d}) \rightarrow \R$ and $h \in {\cal S}(\R^{d})$ we define the derivative that will be denoted either by $D_{h} Y$ or $\partial_{h} Y$ :
\begin{equation}
\label{Derivative}
D_{h} Y = \partial_{h}Y(\varphi)= \lim_{\varepsilon \rightarrow 0} \frac{1}{\varepsilon} (Y(\varphi+h) - Y(\varphi))
\end{equation}
when this limit exists. In the Malliavin calculus terminology this gradient is sometimes called the Malliavin derivative and is shown to be closable with a definition domain $\D_{1}$.
\\
To every action functional ${\cal A}: {\cal S}'(\R^{d}) \rightarrow \R$ we associate the measure $\mu_{{\cal A}}$ defined by :
\begin{equation}
\label{InteractMeasure}
\frac{d\mu_{{\cal A}}}{d\mu_{0}} = R_{{\cal A}}:= \frac{e^{-{\cal A}}}{Z_{{\cal A}}}, \; {\rm with} \; Z_{{\cal A}}= \E e^{-{\cal A}}
\end{equation}
and the expectation with respect to $\mu_{\cal A}$ will be denoted by $\E_{{\cal A}}$. Throughout this section we shall be interested in interacting fields given by a measure of the type (\ref{InteractMeasure}) and the obtained results will be applied to scalar field interactions like $\varphi^{4}$, when we use the regularized measure (\ref{Not8}).
\begin{prop}
\label{IP}
The following integration by parts (I.P.) formulas are valid for every $h\in {\cal S}$ and every functional $Y: {\cal S}'(\R^{d}) \rightarrow \R$ such that $D_{h} Y$ exists $\L^{2}(\P)$:
\begin{equation}
\label{IP1}
 \E D_{h} Y= \E \varphi(C^{-1}h) Y
\end{equation}
\begin{equation}
\label{IP2}
 {\E}_{{\cal A}} D_{h} Y= {\E}_{{\cal A}} \varphi(C^{-1}h) Y + \E_{{\cal A}} Y D_{h}{\cal A}
\end{equation}
They correspond to the free and interacting measure respectively.
\end{prop}
In \cite{GJ87} one can find the following forms of the integration by parts :
\begin{equation}
\label{IP3}
 \E D_{Ch} Y= \E \varphi(h) Y
\end{equation}
\begin{equation}
\label{IP4}
\E \varphi(h) Y= \E (\frac{\delta Y}{\delta \varphi}, Ch) -\E (Y \frac{\delta {\cal A}}{\delta \varphi}, Ch)
\end{equation}
Where the last formula uses the strong differential :
\begin{equation}
\label{IP31}
D_{\delta_{x}} Y= \frac{\delta Y}{\delta \varphi (x)}:= \lim_{\varepsilon \rightarrow 0} \frac{1}{\varepsilon} (Y(\varphi+h) - Y(\varphi))
\end{equation}
when the limit exists and the notation :
\[ ( C\frac{\delta Y}{\delta \varphi}, h)= \int h(y) C(x,y) \frac{\delta Y}{\delta \varphi (y)} dx dy \]
The integration by parts (IP) formula for the free measure (\ref{IP1}) is similar to (\ref{IP3}) while the  IP formula w.r.t. the interacting measure (\ref{IP2}) is slightly different from (\ref{IP4}), in particular because in (\ref{IP2}) we are interested in the expectation w.r.t. the measure $\mu_{{\cal A}}$.
\\
\noindent
{\it Proof.} The quasi-invariance formula (\ref{Quasiinvariance1}) implies that:
\begin{equation}
\label{IP11}
\int Y d\mu(\varphi) Y= \int Y(\varphi+t h) e^{[t(\varphi, C^{-1}h) + t^{2}(h, C^{-1}h)/2]}
\end{equation}
Then (\ref{IP1}) follows by taking the derivative of (\ref{IP11}) at $t=0$. The The quasi-invariance formula (\ref{Quasiinvariance1}) also implies that:
\[  \E_{{\cal A}} Y= \int Y R_{{\cal A}} d\mu(\varphi)) 
=  \int Y (\varphi+t h) \frac{e^{-{\cal A} (\varphi+t h)} }{Z_{{\cal A}}} e^{[t(\varphi, C^{-1}h) + t^{2}(h, C^{-1}h)/2]} d\mu(\varphi) \]
Taking the derivative at $t=0$, this yields :
\[ \E [ D_{h} Y - Y D_{h} {\cal A} - (\varphi, C^{-1}h) Y ] \frac{e^{-{\cal A}} }{Z_{{\cal A}}}=0   \]
which is the formula (\ref{IP2}). $\Box$
\\
\\
\noindent
$\bullet$ {\bf Notations and examples :}
\\
\noindent
We shall use the following notations : for $n\geq 1, f, h \in {\cal S}$:
\begin{equation}
h(f)= \int h f = \int_{\R^{d}} f(x) h(x) dx, \; \; h_{n}(x) = h(\delta_{x}^{n})= \int_{\R^{d}} \delta_{x}^{n}(y) h(y) dy
\end{equation}
For convenience we shall use the symbol $\partial_{h}$ instead of $D_{h}$, that is, we set :
\[ \partial_{h} Y = D_{h} Y,\]
and there will be no confusion with the usual partial derivative $\partial_{x_{i}}$, as the $x_{i}\in \R$ and $h\in {\cal S}$

For a covariance operator $C$ (typically $C=(-\Delta+1)^{-1}$, but the results are valid for other cases) we set :
\begin{equation}
\f= C f, \; \; \sigma_{f}= \E \varphi(f)^{2}=(f,\f)=(f, Cf)
\end{equation}
When the action $\A$ is $\A_{n}$ we set: $ \E_{n}:=\E_{\A}$. The basic formula (\ref{IP2}) will be:
\begin{equation}
\label{MC1}
 \E_{n}\ph1  Y= \E_{n} \varphi(C^{-1}h) Y + \E_{n} Y \ph1 \A_{n}
\end{equation}
and by differentiating w.r.t. $\h$ instead of $h$:
\begin{equation}
\label{MC2}
 \E_{n}\partial_{\h}  Y= \E_{n} \varphi(h) Y + \E_{n} Y \partial_{\h} \A_{n}
\end{equation}
These two formulas and the following examples will be repeatedly used in this paper.
\\
\\
\noindent
$\bullet$ {\it Examples: }
\noindent
\begin{equation}
\label{MC3}
Y(\varphi)= \varphi(f)^{p} \rightarrow \ph1 \varphi(f)^{p}= p \times h(f) \varphi(f)^{p-1}, \; \text{in particular:}\; \ph1 \varphi(f)= h(f).
\end{equation}
With $f=\delta_{x}^{n}$:
\begin{equation}
\label{MC4}
Y(\varphi)= \varphi(\delta_{x}^{n})^{p}=\varphi_{n}(x)^{p} \rightarrow \ph1 \varphi_{n}(x)^{p}= p h(\delta_{x}^{n}) \varphi_{n}(x)^{p-1}=p h_{n}(x) \varphi_{n}(x)^{p-1}
\end{equation}
We can also verify that:
\begin{equation}
\label{MC5}
Y(\varphi)= :\varphi(f)^{p}: \rightarrow \ph1 :\varphi(f)^{p}:= p h(f) :\varphi(f)^{p-1}:,
\end{equation}
and when we deal with the gradient $:(\nabla \varphi(x))^{2}:$, we have:
\begin{equation}
\label{MC6}
\ph1 (\partial_{x_{i}} \varphi_{n}(x))^{2}= 2 \partial_{x_{i}} \varphi_{n}(x) \partial_{x_{i}} (\partial_{h} \varphi_{n}(x))
                                     =2 \partial_{x_{i}} \varphi_{n}(x) \partial_{x_{i}} h_{n}(x)
\end{equation}
and
\begin{equation}
\label{MC7}
\int \ph1 (\partial_{x_{i}} \varphi_{n}(x))^{2}= \int 2 \partial_{x_{i}} \varphi_{n}(x) \partial_{x_{i}} h_{n}(x) = 2 \partial_{x_{i}} \varphi_{n}(\partial_{x_{i}} h_{n})
\end{equation}
\begin{prop}
\label{PropMC1}
For all $f, h\in {\cal S}$, we have the following dynamic equations for the 2d moment and the two-points function:
\begin{equation}
\label{MC8}
\E_{\A} \varphi(f)^{2}= \sf - \E_{\A} \varphi(f) \pr \A
\end{equation}
\begin{equation}
\label{MC8b}
\E_{\A} \varphi(f) \varphi(h)= (\f,h) - \E_{\A} \varphi(f) \partial_{h} \A
\end{equation}
We also have for the higher order moments:
\begin{equation}
\label{MC9}
\E_{\A} \varphi(f)^{p}= (p-1) \sf \E_{\A}\varphi(f)^{p-2} - \E_{\A} \varphi(f)^{p-1} \pr \A
\end{equation}
and the following expression for the 2d moments:
\begin{equation}
\label{MC10}
\E_{\A} \varphi(f)^{2}= \sf + \E_{\A} ((\pr \A)^{2}- \pr^{2} \A)
\end{equation}
\end{prop}
{\it Proof.} We consider the random variable $Y=\varphi(f)$ and we have :
\begin{equation}
\E_{\A} \pr \varphi(f) =
\begin{cases}
(f,\f)=\sf & \text{by differentiation}\\
\E_{\A} \varphi(f) \varphi(C^{-1}\f) + \E_{\A} \varphi(f) \pr \A  & \text{by integration by parts (\ref{IP2})}
\end{cases}
\end{equation}
This gives (\ref{MC8}) because $C^{-1}\f=C^{-1} C f=f$. The second formula (\ref{MC8b}) is obtained in the same way by taking $Y=\varphi(h)$.
The formula (\ref{MC9}) of the higher order moments is also obtained by taking $Y=\varphi(f)^{p-1}$.

The expression (\ref{MC10}) can also be obtained by the above calculus : we take $Y=\pr \A$ and we have:
\begin{equation}
\E_{\A} \pr \A =
\begin{cases}
(f,\f)=\pr^{2} \A & \text{by differentiation}\\
\E_{\A} \varphi(C^{-1}\f) \pr \A + \E_{\A} \pr \A \pr \A  & \text{by the I.P. formula (\ref{IP2})}
\end{cases}
\end{equation}
which gives (\ref{MC10}) since $\E_{\A} \varphi(C^{-1}\f) \pr \A = \E_{\A} \varphi(f) \pr \A= \sf - \E_{\A}\varphi(f)^{2}$ by (\ref{MC8}). $\Box$

\subsection{Generating functional and moments of the interacting field}

We define the generating functional of the regularized interacting field measure $\mu_{{\cal A}}$ by :

\begin{equation}
\label{GF1}
G_{\A}(t,f) = \E_{\A} e^{t\varphi(f)}
\end{equation}
The corresponding moments of the fields are  
\begin{equation}
\label{GF2}
S_{\A,p}(f) := \E_{\A} \varphi(f)^{p} = \frac{d}{dt} G_{\A}(t,f)\vert_{t=0}
\end{equation}
Using the quasi-invariance property we get:

\begin{equation}
\label{GF3}
G_{\A}(t,f) = \frac{1}{Z_{\A}}\E e^{t\varphi(f) + t\h(f)} e^{-{\cal A}(\varphi +th)} e^{[t(\varphi, C^{-1}h) + t^{2}(h, C^{-1}h)/2]},
\end{equation}
and provided we choose $h=\f=C f $ we will have:
\begin{equation}
\label{GF4}
G_{\A}(t,f) = \E_{n} e^{U_{\A}(t,f)}:= \E_{\A} e^{-{\cal A}(\varphi+t\f)+t^{2} (f,\f)/2}
\end{equation}
We set:
\begin{equation}
\label{GF4}
U_{\A}(t,f):= {\cal A}(\varphi+t\f) -\sigma_{f} \frac{t^{2}}{2}, \; \; {\rm with: } \; \sigma_{f}=(f,\f)=(f,Cf)
\end{equation}
We have in particular :
\begin{align}
\label{GF5}
\partial_{0} U_{\A}(t,f) := & \frac{\partial }{\partial t}\vert_{t=0}  U_{\A}(t,f) = \pr \A_{n} \\
\partial_{0}^{2} U_{\A}(t,f) := & \frac{\partial^{2} }{\partial t^{2}}\vert_{t=0}  U_{\A}(t,f) = -\sigma_{f} + \pr \A^{2} \\
\partial_{0}^{p} U_{\A}(t,f) := & \frac{\partial^{p} }{\partial t^{p}}\vert_{t=0}  U_{V}(t,f) = \pr^{p} \A, \; p\geq 3
\end{align}
\noindent
$\bullet$ {\bf First moments of the regularized interacting field}
\\
\noindent
These moments are given by :
\begin{equation}
\label{GF6}
\E_{\A} \varphi(f)^{p}= \frac{\partial^{p} }{\partial t^{p}}\vert_{t=0} G_{\A}(t,f)
\end{equation}
Then we have for the second moments :
\begin{align}
\label{Moment2}
\E_{\A} \varphi(f)^{2} = & \E_{\A}(- \partial_{0}^{2} U_{\A}(t,f) + (\partial_{0} U_{\A}(t,f))^{2} \nonumber \\
               = & \sigma_{f} + \E_{\A}(-\pr^{2} \A + (\pr \A)^{2})
\end{align}
and the 4th moments :
\begin{align}
\label{Moment4}
\E_{\A} \varphi(f)^{4} = & \E_{\A}\{- \p^{4} U_{\A} + 3 (\p^{2} U_{\A})^{2} + 4\p^{3} U_{\A}\p U_{\A} -6 \p^{2} U_{\A}(\p U_{\A})^{2} + (\p U_{\A})^{4} \} \nonumber  \\
               = & \E_{\A}\{- \p^{4} \A + 3 (\p^{2} \A)^{2} + 4\p^{3} \A\p \A -6 \p^{2} \A(\p \A)^{2} + (\p \A)^{4} \} \nonumber  \\
							+ & 3 \sigma_{f}^{2} + 6\sf ((\p \A)^{2} - (\p^{2} \A))  \}
\end{align}
where we have omitted the arguments $(t,f$ of the functional $U_{\A}$ and the subscript $0$ of $\p_{0}U$. The 6th moments are given by:
\begin{align}
\label{Moment6}
\E_{\A} \varphi(f)^{6} = & \E_{\A}\{ (\p U_{\A})^{6}+15 \p^{2} U_{\A} \p^{4} U_{\A}- (\p U_{\A})^{2} \p^{4} U_{\A} + 10 (\p^{3} U_{\A})^{2} - 15 (\p^{2} U_{\A})^{3} \nonumber  \\
               - & 60 \p U_{\A} \p^{2} U_{\A} \p^{3} U_{\A} + 45 (\p U_{\A})^{2}(\p^{2} U_{\A})^{2} + 20 (\p U_{\A})^{3}\p^{3} U_{\A} \nonumber  \\
							 - & 15 (\p U_{\A})^{4} \p^{2} U_{\A}, \}
\end{align}
and using (\ref{GF5}) :
\begin{align}
\label{Moment61}
\E_{\A} \varphi(f)^{6} = & \E_{\A}\{ (\p \A)^{6}+15 \p^{2} \A \p^{4} \A- (\p \A)^{2} \p^{4} \A + 10 (\p^{3} \A)^{2} 
                        - 15 (\p^{2} \A_{n})^{3} \nonumber  \\
               - & 60 \p \A \p^{2} \A \p^{3} \A + 45 (\p \A)^{2}(\p^{2} \A)^{2} + 20 (\p \A)^{3}\p^{3} \A  \nonumber  \\
							 - & 15 (\p A)^{4} \p^{2} \A, \} + {\rm complementary \;  terms}
\end{align}
where:
\begin{align}
\label{Moment62}
{\rm complementary \; terms} = & \E_{\A}\{ 15 \sf((\p \A)^{4} - \p^{4} \A) + 15 (3 \sf (\p^{2} \A)^{2} + 3 \sf^{2} \p^{2} \A \nonumber \\
                           - & \sf^{3} + 60 \sf \p \A \p^{3} \A + 45 (\sf^{2}- 2 \sf \p^{2} \A) (\p \A)^{2} \}
\end{align}

\section{Lower bound revisited and the least action principle for QFT}
In this section we shall study some refinements of the lower bound of large deviations of the action $\A_{n}$ obtained in \cite{AA1}, which will lead to the formulation of the principle of least action in Euclidean field theory. The notations and conventions of \cite{AA1} will be used. We start from the action :
\begin{align}
\label{LB1}
{\cal A}_{n}= & g_{n} \int_{V}: \varphi_{n}^{4}:+ m_{n} \int_{V}: \varphi_{n}^{2}:+a_{n} \int_{V}: (\partial\varphi_{n})^{2}:\\
            = & g_{n} c_{n}^{2} \int_{V}: \psi_{n}^{4}:+ m_{n}c_{n} \int_{V}: \psi_{n}^{2}: 
						       + a_{n} c_{n} d_{n} \frac{1}{d_{n}} \int_{V}: (\partial\psi_{n})^{2}:  \nonumber
\end{align}
where 
\begin{equation}
\label{LB1b}
\psi_{n} (x) = \varphi_{n}(x) / (E\varphi_{n}(x)^{2})^{1/2} = \varphi_{n}(x) / c_{n}^{1/2}, \; \;  d_{n} = \E (\partial\psi_{n})^{2}(0) \sim K n^{d-2}
\end{equation}
\noindent
We set :  
\begin{equation}
\label{LB2}
I_{n}= \int_{V}: \psi_{n}^{4}:, \; \;     M_{n}=\int_{V}: \psi_{n}^{2}:, \; \;   D_{n}=\frac{1}{d_{n}}\int_{V}: (\partial\psi_{n})^{2}:
\end{equation}
\noindent
We shall use the following convention for the finite volume: we take $V=[0,1]^{d}$ and for $n\geq 1, i=(i_{1}, ..., i_{d}), i_{k}=0, ..., n-1$, we set:
 $$V_{ni}= [i_{1}, i_{1}+1]\times [i_{2}, i_{2}+1] \times ...[i_{d}, i_{d}+1] $$
For definiteness we limit ourselves to this configuration, but other types of volumes could be considered as well.
\begin{Rk}
We recall the role of the factor $d_{n}$: first we have $\E (\partial\varphi_{n}(x))^{2}\sim K_{1} n^{d} \sim K_{2} n^{2}c_{n}$, with $c_{n}= \E \varphi_{n}(x)^{2} \sim K n^{d-2}$. We also have $\E (\partial\psi_{n}(x))^{2}= K_{1} n^{2}$ and therefore $:(\partial\psi_{n}(x))^{2}: = (\partial\psi_{n}(x))^{2} - K_{1} n^{2} $, and when we write the 3d term of the Lagrangian as a mean of series, i.e.:
\begin{equation}
\label{LB2b}
\frac{1}{d_{n}}\int_{V}: (\partial\psi_{n})^{2}: = \frac{1}{n_{d}}\sum_{i=1}^{n^{d}} \tilde{D}_{n,i}
\end{equation}
we are led to take : $\tilde{D}_{n,i} = (1/d_{n}) \int_{V_{ni}}: (\partial\psi_{n})^{2}(\frac{x}{n}):dx $ in order to have $  0< C_{1} \leq \E \tilde{D}_{n,i}^{2} \leq C_{2} <+\infty $.
\end{Rk}

The outcome depends on the dominant terms of the action and this leads to consider the following cases and to define a secondary action ${\cal A}_{n}^{(1)}$:

$\bullet$ {\bf Case 1:} $g_{n}c_{n}^{2} \geq K m_{n}c_{n}$ and $g_{n}c_{n}^{2} \geq K' a_{n}c_{n}d_{n}$ for some $K, K'$. In this case we write :

\begin{align}
\label{LB3}
{\cal A}_{n}(\varphi_{n})= & g_{n} c_{n}^{2} [\int_{V}: \psi_{n}^{4}:+ \frac{m_{n}}{g_{n} c_{n}} \int_{V}: \psi_{n}^{2}: 
						       + \frac{a_{n}d_{n}}{g_{n} c_{n}} \frac{1}{d_{n}} \int_{V}: (\partial\psi_{n})^{2}: ]\\
            = & g_{n} c_{n}^{2} [ \lambda_{n} \int_{V}: \psi_{n}^{4}:+ \alpha_{n} \int_{V}: \psi_{n}^{2}: 
						       + \beta_{n}\frac{1}{d_{n}} \int_{V}: (\partial\psi_{n})^{2}: ] \nonumber \\
						= & g_{n} c_{n}^{2} {\cal A}_{n}^{(1)} (\psi_{n})  \nonumber
\end{align}

with $\lambda_{n}=1, \alpha_{n} = m_{n}/ g_{n}c_{n}$ and $\beta_{n} = a_{n}d_{n}/ g_{n}c_{n}$. 

$\bullet$ {\bf Case 2:} $ m_{n} \geq K g_{n}c_{n}$ and $m_{n} \geq K' a_{n}d_{n}$ for some $K, K'$. In this case we write :
\begin{align}
\label{LB4}
{\cal A}_{n} (\varphi_{n})= & m_{n}c_{n}  [ \lambda_{n} \int_{V}: \psi_{n}^{4}:+ \alpha_{n} \int_{V}: \psi_{n}^{2}: 
						       + \beta_{n}\frac{1}{d_{n}} \int_{V}: (\partial\psi_{n})^{2}: ]  \\
						= & m_{n}c_{n}  {\cal A}_{n}^{(1)} (\psi_{n}) \nonumber
\end{align}
with $\lambda_{n}= g_{n} c_{n}/ m_{n}, \alpha_{n} = 1$ and $\beta_{n} = a_{n}d_{n}/ m_{n}$.

$\bullet$ {\bf Case 3:} $ a_{n}d_{n} \geq K g_{n}c_{n}$ and $a_{n}d_{n} \geq K' m_{n}$ for some $K, K'$. In this case we write :
\begin{align}
\label{LB5}
{\cal A}_{n} (\varphi_{n})= & a_{n}c_{n}d_{n}  [ \lambda_{n} \int_{V}: \psi_{n}^{4}:+ \alpha_{n} \int_{V}: \psi_{n}^{2}: 
						       + \beta_{n}\frac{1}{d_{n}} \int_{V}: (\partial\psi_{n})^{2}: ] \\
						= a_{n}c_{n}d_{n}  {\cal A}_{n}^{(1)}(\psi_{n})  \nonumber
\end{align}
with $\lambda_{n}= g_{n} c_{n}/ a_{n}d_{n}, \alpha_{n} = m_{n}/ a_{n}d_{n}$ and $\beta_{n} = 1$.

Let $\Lambda^{*}(x)=\limsup_{n} \Lambda^{*}_{n}(x)$, and $\Lambda^{*}_{n}$ be the Legendre-Fenchel transform of the function $\Lambda_{n}$ defined by $\Lambda_{n}(\theta)= \log \E \exp(-\theta X_{n,i})$, with :
\begin{equation}
\label{LB6}
X_{n,i} = \lambda_{n} I_{n,i}+\alpha_{n} M_{n,i}+\beta_{n} D_{n,i}
\end{equation}
where :
\[ I_{n,i}= \int_{V_{ni}}: \psi_{n}^{4} (\frac{x}{n}):dx ,\;\; M_{n,i}=\int_{V_{ni}}: \psi_{n}^{2}(\frac{x}{n}):dx,\;\; 
D_{n,i}=\int_{V_{ni}}: (\partial\psi_{n})^{2}(\frac{x}{n}):dx\]
We also write :
\begin{align}
\label{LB6bis}
X_{n,i} & = \int_{V_{ni}}: l_{n}(\psi_{n}(\frac{x}{n}))dx, \\ \nonumber 
l_{n}(x) &= \lambda_{n} : \psi_{n}^{4}(x):+\alpha_{n} :\psi_{n}^{2}(x): +\beta_{n}: (\partial\psi_{n})^{2}(x).
\end{align}
\\
\noindent
$\bullet$ {\bf From a large deviation lower bound to the principle of least action}

The main result of \cite{AA1} is that for some rate function $\Lambda^{*}$ with non-empty domain we have:
\begin{equation}
\label{LB7}
\liminf_{n\rightarrow \infty} \frac{1}{n^{d}}\log {\rm Pr} (-\A_{n}^{(1)} (\psi_{n}) \in [a,b]) \geq  -\inf_{x\in ]a,b[} \Lambda^{*}(x)
\end{equation}
Let $z^{max} > 0$ be such that:
\begin{equation}
\label{LB7b}
{\rm Domain} (\Lambda^{*}) \bigcap \R^{+}=[0, z^{max}[
\end{equation}
For any $\varepsilon > 0$, there is a $n_{\varepsilon}$ such that for $n\geq n_{\varepsilon}$:
\begin{equation}
\label{LB8}
 {\rm Pr} (-\A_{n}^{(1)}(\psi_{n}) \in [a,b]) \geq  e^{-n^{d}(\inf_{x\in ]a,b[} \Lambda^{*}(x)+\varepsilon)}
\end{equation}
In particular for $[a,b[=[z,\infty[$
\begin{equation}
\label{LB9}
 {\rm Pr} (-\A_{n}^{(1)} (\psi_{n}) \geq z) \geq  e^{-n^{d}(\inf_{x>z} \Lambda^{*}(x)+\varepsilon)}
\end{equation}
 (\ref{LB8}) is of course interesting only when $z\in {\rm Domain} (\Lambda^{*})$ for if $\Lambda^{*}(z)=\infty$, (\ref{LB4}) means that $ {\rm Pr} (-\A_{n}^{(1)}(\psi_{n}) \geq z) \geq 0$ which is trivial. The random variables $X_{n,i}$ have the same law for each $n$;  let $Y_{n}$ be a random variable having the law of $-X_{n,i}$ and :
\[ \Lambda_{n}(\theta)= \log \E \exp(\theta Y_{n}), \; \; \Lambda_{n}'(\theta)= \frac{\E Y_{n}e^{\theta XY_{n}}}{\E e^{\theta Y_{n}}} \]
In our case $\Lambda_{n}, \Lambda^{*}_{n}$ are continuous and, $\Lambda_{n}$ being convex, we have : $\Lambda_{n}'$ is increasing and continuous. By the gneral properties of the Legendre-Fenchel transform we also have : $\Lambda_{n}'(\theta)=y \Rightarrow \Lambda^{*}_{n}(y)= \theta y - \Lambda_{n}(\theta) $. 
Hence $\Lambda^{*}_{n}$ is defined on $[0, \Lambda^{ '}_{n}(\infty)[$, with $\Lambda^{'}_{n}(\infty)=\lim_{\theta\rightarrow \infty} \Lambda^{ '}_{n}(\theta)$.
Now what would be the value of $\lim_{\theta\rightarrow \infty} \Lambda^{* '}_{n}(\theta)$ ? 

\noindent
The fact that will be seen below (Proposition \ref{PropPLA1}) is that :
\begin{equation}
\label{LB10}
\lim_{\theta\rightarrow \infty} \Lambda^{* '}_{n}(\theta)=ess \max (Y_{n}) = ess \max (-X_{n,i}) = z^{max}_{n} 
\end{equation}
where $ess \max $ is the essential maximum of a random variable (defined in Lemma \ref{lemma-MaxRV}). We will also see that this essential maximum $z^{max}_{n}$ of the random variable $-X_{n,i}$ given by (\ref{LB6bis}) corresponds to the maximum of the action:
\begin{equation}
\label{LB11}
- \A^{(1)}_{n}(\chi)=  - \int_{V} [ \lambda_{n}\chi^{4}+ \alpha_{n}\chi^{2}+\frac{\beta_{n}}{d_{n}} (\partial\chi)^{2}]
\end{equation}
Now, to get the principle of the least action, we note that the lower bound (\ref{LB8}) implies that for any $\delta > 0$, $\varepsilon > 0$, there is some $n_{\varepsilon}\geq 1$ such that for all $n \geq n_{\varepsilon}$ :
\begin{equation}
\label{LB12}
 {\rm Pr} (-\A_{n}^{(1)} (\psi_{n}) \geq z^{max}-\delta) \geq  \exp [-n^{d}(\inf_{ \{x> z^{max}-\delta \}} \Lambda^{*}(x)+\varepsilon) ]
\end{equation}
In the case 1 for instance, this yields :
\begin{align}
\label{LB13}
 Z_{n}= E_{n} e^{-\A_{n}(\varphi)} & = E_{n} e^{-g_{n}c_{n}^{2}\A_{n}^{(1)} (\psi_{n})} \\ \nonumber
                                & \geq e^{-g_{n}c_{n}^{2} (z^{max}-\delta)} {\rm Pr} (-\A_{n}^{(1)} (\psi_{n}) \geq z^{max}-\delta)     \\ \nonumber
																& \geq e^{-g_{n}c_{n}^{2} (z^{max}-\delta)} e^{-n^{d}(\inf_{ \{ x> z^{max}-\delta \}} \Lambda^{*}(x)+\varepsilon)} 
\end{align}
and the interaction field measure $\mu_{n}$ will satisfy:
\begin{align}
\label{LB14}
 \frac{d\mu_{n}}{d\mu_{0}}= \frac{e^{-\A_{n}(\varphi)}}{ Z_{n}}& = \frac{e^{-g_{n}c_{n}^{2}\A_{n}^{(1)} (\psi_{n})}}{ Z_{n}}  \\ \nonumber 
                          & \leq \exp \{ g_{n}c_{n}^{2}[ - \A_{n}^{(1)}(\psi) - (z^{max}-\delta)] + n^{d}[\inf_{x> z^{max}-\delta} \Lambda^{*}(x)+\varepsilon] \}  
\end{align}
This last inequality shows that the field measure $\mu_{n}$ is supported in the set where $- \A_{n}^{(1)} (\psi_{n}) \geq (z^{max}_{n}-\delta) $ which means that the the action $- \A_{n}^{(1)}$ is close to the maximum $z^{max}_{n}$. In other words, the action $\A_{n}^{(1)}$ is close to the minimum $-z^{max}_{n}$

Now we shall make precise these ideas and state the principle of the least action. For the $\varphi^{4}$ regularized action (\ref{LB1}) we consider the reduced action defined above :
\begin{equation}
\label{LA1}
\A^{(1)}_{n}(\psi_{n})=  \int_{V}  [ \lambda_{n}:\psi_{n}^{4}:+ \alpha_{n} : \psi_{n}^{2}:+ \frac{\beta_{n}}{d_{n}} : (\partial\psi_{n})^{2}:]
\end{equation}
We note that the 'constant' term or vacuum renormalization in the Lagrangian $\A_{n}$ (\ref{LB1}) is:
\begin{equation}
\label{LA2}
Vac(n)=  (3 g_{n}c_{n}^{2} - m_{n}c_{n} - a_{n} n^{d}) V
\end{equation}
and for the reduced action we have :
\begin{equation}
\label{LA3}
\A^{(1)}_{n}(\psi_{n})=  \int_{V} [ \lambda_{n}\psi_{n}^{4}+ \alpha_{n}\psi_{n}^{2}+ \frac{\beta_{n}}{d_{n}} (\partial\psi_{n})^{2}] + Vac^{(1)}(n)
\end{equation}
with $Vac^{(1)}(n)= (3 \lambda_{n} - \alpha_{n} - \beta_{n} n^{2}/d_{n}) V = (3 \lambda_{n} - \alpha_{n} - \beta_{n}) V $ is a vacuum renormalization constant up to a multiplicative sequence ($g_{n}c_{n}^{2}$ in the case 1, etc.). Then we consider the classical action:
\begin{equation}
\label{LA4}
\A^{c}_{n}(\chi)=  \int_{V} [ \lambda_{n}\chi^{4}+ \alpha_{n}\chi^{2}+\frac{\beta_{n}}{d_{n}} (\partial\chi)^{2}].
\end{equation}
Since the sequences $\lambda_{n}, \alpha_{n}, \beta_{n}$ are bounded by construction in the cases 1-3, we may suppose that they have limits  $\lambda, \alpha, \beta$ and we consider the following classical action:
\begin{equation}
\label{LA5}
\A^{c,2}_{n}(\chi)=  \int_{V} [ \lambda\chi^{4}+ \alpha \chi^{2}+\frac{\beta}{d_{n}} (\partial\chi)^{2}]
\end{equation}

\begin{thm}
\label{ThmPLA}
The $\varphi^{4}$ quantum field model is considered in the Euclidean framework through the regularized action in a finite volume :
\begin{equation}
\label{LA6}
{\cal A}_{n}= g_{n} \int_{V}: \varphi_{n}^{4}:+ m_{n} \int_{V}: \varphi_{n}^{2}:+a_{n} \int_{V}: (\partial\varphi_{n})^{2}:
\end{equation}
We add further the following :
\\
\noindent
Assumption (A) 
\[ \; g_{n}c_{n}^{2}/n^{d}\sim g_{n}n^{d-4}  \rightarrow \infty \; {\rm or } 
    \; m_{n}c_{n}/n^{d} \sim m_{n}/n^{2}\rightarrow \infty \; {\rm or }
		\; a_{n}c_{n} n^{2}/n^{d} \sim  a_{n} \rightarrow  \infty  \]
Then, for large $n$, the measure $\mu_{n}$ associated with the regularized $\varphi^{4}$ model est supported in the sets :
\begin{equation}
\label{LA7}
\Sigma_{\varepsilon,n}^{1} = \{ \varphi \in {\cal S}'(\R^{d}): {\cal A}^{c}_{n} (\psi_{n}) \leq \min_{\chi \in C^{2}(\R^{d})}{\cal A}_{n}^{c} (\chi) + \varepsilon, \; {\rm with}: \; \psi_{n}=\varphi_{n}/\sqrt{c_{n}} \}   
 \end{equation}
\begin{equation}
\label{LA8}
\Sigma_{\varepsilon,n}^{2} = \{ \varphi \in {\cal S}'(\R^{d}): {\cal A}^{c,2}_{n} (\psi_{n}) \leq \min_{\chi \in C^{2}(\R^{d})}{\cal A}_{n}^{c,2} (\chi) + \varepsilon, \; {\rm with}: \; \psi_{n}=\varphi_{n}/\sqrt{c_{n}} \}   
 \end{equation}

That the measure $\mu_{n}$ is supported in the previous sets means that for all $\varepsilon$, there is some $n_{1}\geq 1$ such that the probabilities $Pr(\varphi \in \Sigma_{\varepsilon}^{i, c}), i=1,2$ is exponentially small, i.e. $\leq \exp(-nk(\varepsilon))$ for some function $k(\varepsilon)$. 
\end{thm} 
Since we already have ${\cal A}^{c}_{n} (\psi_{n}) \geq \min_{\chi \in C^{2}(\R^{d})}{\cal A}_{n}^{c} (\chi)$, this theorem means that the normalized limiting field (more precisely the field $\varphi/\sqrt{c_{n}}$ under the measure $\mu_{n}$, with $n$ being large) is actually supported around the minimum of the action ${\cal A}_{n}^{c} (\chi)$.
\\
\noindent
The previous theorem provides a rigorous link between quantum field theory and classical field theory, and can be extended to more general scalar fields ($P(\varphi)_{d}$, etc.). This will be discussed elsewhere in more detail.  
\\
\noindent
Equivalently, the measure $\mu_{n}$ associated with the regularized $\varphi^{4}$ model est supported in the set:
\begin{equation}
\label{LA9}
 \Sigma_{\varepsilon,n}^{3} = \{\varphi \in {\cal S}'(\R^{d}): {\cal A}^{1}_{n} (\psi_{n}) \leq \min_{\chi \in C^{2}(\R^{d})} {\cal A}^{1}_{n}(\chi) + \varepsilon, \; {\rm with}: \; \psi_{n}=\varphi_{n}/\sqrt{c_{n}}  \}
\end{equation}

\begin{Rk}
The assumption (A) is always fulfilled if $d\geq 5$, unless the coupling constant is vanishing $g_{n}\rightarrow 0$, in which case the limiting field is Gaussian (by the skeleton inequalities if we work with the lattice approximation). In dimension 4, this assumption means that we should have either $g_{n}\rightarrow \infty$ or $m_{n}/n^{2}\rightarrow \infty$ or $a_{n}\rightarrow \infty$, a condition that can presumably be relaxed.
\end{Rk}
{\it Proof.} 
We start from the lower bound (\ref{LB7}): let $[0,z^{max}[ = {\rm Domain}(\Lambda^{*})\bigcap \R^{+}$; we may have $z^{max}=+\infty$. Let $\delta > 0$ be a fixed positive number. Then, for any $\varepsilon > 0$ fixed, there is some $n_{\varepsilon} \geq 1$ such that:
\begin{equation}
\label{LA20}
 {\rm Pr} (-\A^{(1)}_{n} (\psi_{n}) \geq r_{\infty} - \delta) \geq  \exp [-n^{d}(\sup_{\{x>r_{\infty}-\delta\}} \Lambda^{*}(x)+\varepsilon)]
\end{equation}
(and we have $\sup_{x>r_{\infty}-\delta} \Lambda^{*}(x) < +\infty$). If we show that $z^{max} - \delta$ is close to $y^{min}$ up to any small $\delta'> 0$, where  $y^{min}_{n}$ is the minimum of the action $\A^{c}_{n}$ then the theorem will be proved according to the ideas outlined before its statement (\ref{LB7})-(\ref{LB14}). To do this we proceed as follows: Let $y_{n}= \Lambda^{'}_{n}(\infty):=\lim_{\theta\rightarrow \infty}\Lambda^{'}_{n}(\theta)$ if it exists. Then we have $[0, y_{n}[ \subset Domain (\Lambda^{*}_{n}) $ and we will show that:
\\
\noindent
$\bullet$ (a) For large $n$ and at least for some subsequence (still denoted by $n$), we have :
\begin{equation}
\label{LA21}
 y_{n} - \delta \in Domain (\Lambda^{*})
\end{equation}
\noindent
$\bullet$ (b) $y_{n}$ exists and corresponds to the minimum of the action $\A^{(1)}_{n}$ which is also the minimum of the action $\A^{c}_{n}$.
\\
\\
\noindent
Let $[0, r_{n}[$ be the definition domain of $\Lambda^{*}_{n}$. We have $y_{n}= \Lambda^{'}_{n}(\infty) \leq r_{n}$.
The assertion (a) follows from the fact that $\Lambda^{*}(x) = \limsup_{n} \Lambda^{*}_{n}(x)$, which implies that for any $\varepsilon > 0 $ there exists $n_{2}\geq 1 $ such that : $ \Lambda^{*}_{n}(x) \leq \Lambda^{*}(x) + \varepsilon, \forall n\geq n_{2} $  and $ \Lambda^{*}(x) - \varepsilon \leq \Lambda^{*}_{n}(x) $ for infinitely many $n\geq n_{2}$. The last inequality implies that $\Lambda^{*}(x) $ is bounded ($ \leq  \Lambda^{*}_{n}(x) + \varepsilon$) as soon as $\Lambda^{*}_{n}(x) $ is. In particular for the later infinitely many $n$, since $y_{n} \leq r_{n}$ and thus $y_{n}-\delta \in Domain (\Lambda^{*}_{n})$ we will have $y_{n}-\delta \in Domain (\Lambda^{*})$. This proves (a).

The assertion (b) will be proved in two steps:
\\
\noindent
(b.1) We shall show that:
\begin{equation}
\label{ToBeProved}
{\Lambda'}_{n}(\infty)= \lim_{\theta \rightarrow \infty} \frac{\E Y_{n} e^{\theta Y_{n}}}{\E e^{\theta Y_{n}}} = {\rm ess}\max Y_{n}
\end{equation}
with $Y_{n}$ being a random variable that has the same law as :
\begin{equation}
\label{LA22}
-X_{n,i} = -\int_{\Lambda_{ni}}: l_{n}(\psi_{n}(\frac{x}{n}))dx
\end{equation}
(cf. (\ref{LB6bis}) and ${\rm ess}\max Y_{n}$ is the essential maximum of the random variable $Y_{n}$, see lemma \ref{lemma-MaxRV}. This part is not trivial and is the subject of the proposition \ref{PropPLA1} below. We suppose for the moment that (b.1) is valid.
\\
\noindent
(b.2) In this second step we make the link between the essential maximum of the variable $Y_{n}$ which corresponds to the essential minimum of the r.v. $X_{n,i}$ and the minimum of the actions $\A^{(1)}$ or $\A^{c}$ of the theorem. 

By (a) we have for large $n$ belonging to the subsequence mentioned above $r_{n} - \delta \in Domain (\Lambda^{*})$ and $\Lambda^{'}_{n}(\infty) - \delta \in Domain (\Lambda^{*})$ since $\Lambda^{'}_{n}(\infty) \in [0, r_{n} [$ the definition domain of  $\Lambda^{*}_{n}$ . This implies that for sufficiently large $n$, we have $z^{max} \geq \max_{\omega} Y_{n}(\omega) -\delta $. On the other hand, for a given $n$, $\max_{\omega} Y_{n}(\omega)$, which is the essential maximum of :
\begin{equation}
\label{LA23}
-{\cal A}^{1}_{n} (\psi)= -[\lambda_{n} \int_{V_{ni}}:\psi^{4}_{n}(\frac{x}{n})):dx+ \alpha_{n} \int_{V_{ni}} :\psi^{2}_{n}(\frac{x}{n})):dx
                            +\frac{\beta_{n}}{d_{n}} \int_{V_{ni}}:(\partial\psi_{n})^{2}(\frac{x}{n})):dx]
\end{equation}
when paths $\varphi$ or $\psi$ run over ${\cal S}'(R^{d})$ or the regularized paths $\psi_{n}$ run over ${\cal S}(R^{d})$ is the same as the essential maximum of :
\begin{equation}
\label{LA24}
-{\cal A}^{c}_{n}(\chi)= -[\lambda_{n} \int_{V_{ni}}:\chi^{4}(x)):dx+ \alpha_{n} \int_{V_{ni}} :\chi^{2}(x):dx
                            +\frac{\beta_{n}}{d_{n}}\int_{V_{ni}}:(\partial\chi)^{2}(x):dx]
\end{equation}
where this time the $\chi$ run over the regular path space ${\cal S}(R^{d})$ or $C^{2}(\R^{d})$ and the integration volume can be taken to be $\Lambda$. 

This maximum can be made close up to $\delta$ to the maximum of :
\begin{equation}
\label{LA25}
-{\cal A}^{c}_{n} (\chi)= -[\lambda \int_{\Lambda}\chi^{4}+ \alpha \int_{\Lambda} \chi^{2}+\frac{\beta}{d_{n}} \int_{\Lambda}(\partial\chi)^{2}],
\end{equation}
This is due to the assumption that $\lambda_{n}, \alpha_{n}, \beta_{n}$ converge to $\lambda, \alpha, \beta$ and the fact that the maximum of functionals like ({\ref{LA24}) is continuous w.r.t to the parameters $\lambda, \alpha, \beta$, a fact that is easy to establish. 
This maximum corresponds of course to the minimum of the action ${\cal A}^{c}_{n} (\chi)$.

So let $y^{(m)}_{n} = \max_{\chi} (-{\cal A}^{c}_{n}) = - \min_{\chi} ({\cal A}^{c}_{n}) $. By (a) and (b) $y^{(m)}_{n} - \delta \in {\rm Domain} (\Lambda^{*})$ ;  Then  (\ref{LA20}) means that for all $n\geq n(\delta, \varepsilon)$:
\begin{equation}
\label{LA26}
Pr (-{\cal A}^{1}_{n} (\psi_{n}) \geq y^{(m)}_{n} - 2 \delta) \geq  e^{-n^{d}(\sup_{x>z^{max}-\delta} \Lambda^{*}(x)+\varepsilon)}
\end{equation}
 
$\bullet$ Let us consider the case 1 ($g_{n}c_{n}^{2} \geq K m_{n}c_{n}$ and $g_{n}c_{n}^{2} \geq K' a_{n}c_{n}$ for some $K, K'$). We will have :
\begin{align}
\label{LA27}
Z_{n}= & E e^{g_{n}c_{n}^{2}(-{\cal A}^{1}_{n} (\psi_{n})) } \nonumber \\
    \geq  & = E e^{g_{n}c_{n}^{2}(-{\cal A}^{1}_{n} (\psi_{n})) } 1_{{\cal A}^{1}_{n} (\psi_{n}) \geq y^{(m)}_{n} - 2\delta}  \nonumber \\
		\geq & \exp [g_{n}c_{n}^{2} (y^{(m)} - 2\delta)] \exp [-n^{d}(\sup_{x>z^{max}-\delta} \Lambda^{*}(x)+\varepsilon)] \nonumber \\
		\geq & \exp g_{n}c_{n}^{2} [ y^{(m)} - 2\delta - (\Lambda^{*}(z^{max}-\delta)+\varepsilon) n^{d}/g_{n}c_{n}^{2}  ]
\end{align}
This yields:
\begin{align}
\label{LA28}
 \frac{d\mu_{n}}{d\mu_{0}} & = \frac{e^{-\A_{n}(\varphi)}}{ Z_{n}} = \frac{e^{-g_{n}c_{n}^{2}\A_{n}^{(1)} (\psi_{n})}}{ Z_{n}}  \\ \nonumber
                          & \leq \exp \{ g_{n}c_{n}^{2}[ - \A_{n}^{(1)}(\psi) - (y^{(m)}_{n}-2\delta) + (\Lambda^{*}(z^{max}-\delta)+\varepsilon) n^{d}/g_{n}c_{n}^{2}] \}  \nonumber
\end{align}
This shows that the interaction measure $\mu_{n}$ is supported in the set $\{ - \A_{n}^{(1)}(\psi) \geq (y^{(m)}-2\delta) \}$. For a given random variable $X$ (in $L^{2}({\cal P}), \mu_{n}$ for instance) we have indeed: 
\begin{equation}
\label{LA29}
\E_{n} X = \E_{n} X 1_{\Sigma_{\delta,n}^{1}} + \E_{n} X 1_{\Sigma_{\delta,n}^{1,c}} 
\end{equation} 
($\Sigma_{\delta,n}^{1,c}$ is the complementary set of $\Sigma_{\delta,n}^{1}$) and then, still for $n\geq n(\delta, \varepsilon)$:
\begin{align}
\label{LA30}
(\E_{n} X 1_{\Sigma_{3\delta,n}^{1,c}})^{2}  \leq  &  \E_{n} X^{2}  \E_{n} X  1_{\Sigma_{3\delta,n}^{1,c}}  \\  \nonumber 
		\leq & \E_{n} X^{2}  \exp g_{n}c_{n}^{2} [ y^{(m)}_{n}-2\delta - (y^{(m)}_{n}-\delta) + (\Lambda^{*}(z^{max}-\delta)+\varepsilon) n^{d}/g_{n}c_{n}^{2}] \\ \nonumber 
\end{align}
Now since $\Lambda^{*}(z^{max}-\delta)$ is finite for any given $\delta > 0$ and $n^{d}/g_{n}c_{n}^{2}\rightarrow 0$ by assumption, there exists a $n_{\varepsilon, \delta}^{'} \geq 1$ such that for $n\geq n_{\varepsilon, \delta}^{'}$ the exponential factor in the previous inequality verifies:
\begin{equation}
\label{LA31}
y^{(m)}_{n}-3\delta - (y^{(m)}_{n}-2\delta) + (\Lambda^{*}(z^{max}-\delta)+\varepsilon) n^{d}/g_{n}c_{n}^{2} \leq -\delta/2
\end{equation}
which means:
\begin{equation}
\label{LA32}
(\E_{n} X 1_{\Sigma_{3\delta,n}^{1,c}})^{2}  \leq  \E_{n} X^{2} e^{- g_{n}c_{n}^{2} \delta }
\end{equation}
This proves the theorem in the case 1.
\\
\\
\noindent
$\bullet$ {\bf Case 2:} where $ m_{n} \geq K g_{n}c_{n}$ and $m_{n} \geq K' a_{n}n^{2}$ for some $K, K'$, (\ref{LA27}) becomes:
\begin{equation}
\label{LA33}
Z_{n}=  \E e^{m_{n}c_{n}(-{\cal A}^{1}_{n} (\psi_{n})) } \geq \exp m_{n}c_{n} [ y^{(m)}_{n} - 2\delta - (\Lambda^{*}(z^{max}-\delta)+\varepsilon) n^{d}/m_{n}c_{n}^{2}  ]
\end{equation}
where this time the expression of $(-{\cal A}^{1}_{n})$ is changed as in (\ref{LB3}). The same arguments used in (\ref{LA28})-(\ref{LA32}) give:
\begin{equation}
\label{LA34}
\E_{n} X 1_{\Sigma_{3\delta,n}^{1,c}}  \leq  \E_{n} X^{2} e^{- m_{n}c_{n}\delta }
\end{equation}
provided that $n^{d}/m_{n}c_{n} \rightarrow 0$.
\\
\\
\noindent
$\bullet$ {\bf Case 3:} where $ a_{n}n^{2} \geq K g_{n}c_{n}$ and $a_{n}n^{2} \geq K' m_{n}$ for some $K, K'$, we get in a similar way :
\begin{equation}
\label{LA35}
Z_{n}=  \E e^{a_{n}c_{n}n^{2}(-{\cal A}^{1}_{n} (\psi_{n})) } 
\geq \exp a_{n}c_{n} n^{2}[ y^{(m)}_{n} - 2\delta - (\Lambda^{*}(z^{max}-\delta)+\varepsilon) n^{d}/a_{n}c_{n} n^{2} ]
\end{equation}
and
\begin{equation}
\label{LA36}
\E_{n} X 1_{\Sigma_{3\delta,n}^{1,c}}  \leq  \E_{n} X^{2} e^{a_{n}c_{n}n^{2} \delta }
\end{equation}
provided that $n^{d}/a_{n}c_{n}n^{2} \rightarrow 0$. $\Box$

We turn now to the proof of the assertion (\ref{ToBeProved}). Its meaning is rephrased in the following:

\begin{prop}
\label{PropPLA1}
Let $Y_{n}$ be the random variable :
\begin{equation}
\label{PLA1}
Y_{n}= \int_{V_{ni}} [\alpha_{n} :\psi^{4}_{n}(\frac{x}{n})):+ \beta_{n} :\psi^{2}_{n}(\frac{x}{n})):
                            +\gamma_{n} :(\partial\psi_{n})^{2}(\frac{x}{n})): ]dx
\end{equation}
and suppose that the deterministic functional :
\begin{equation}
\label{PLA2}
B_{n}(\chi)= -[\alpha_{n} \int_{V_{ni}}:\chi^{4}(\frac{x}{n})):dx+ \beta_{n} \int_{V_{ni}} :\chi^{2}(\frac{x}{n})):dx
                            +\gamma_{n} \int_{V_{ni}}:(\partial\chi)^{2}(\frac{x}{n})):dx]
\end{equation}
has a minimum $y^{min}$ attained in a path $\chi^{min}$ (which means that $-B_{n}$ has a maximum $y^{max}=-y^{min}$) then we have :
\begin{equation}
\lim_{\theta \rightarrow \infty} \frac{\E Y_{n} e^{\theta Y_{n}}}{\E e^{\theta Y_{n}}} = y^{max}=-y^{min} = {\rm ess}\max Y_{n}
\end{equation}
\end{prop}
\noindent
This proposition is a consequence of the following lemma:

\begin{lm}
\label{lemma-MaxRV}
Let $Y$ be a real random variable and suppose that it has an essential maximum $y^{max}$ in the following sense:
\begin{equation}
\label{lm1}
Y \leq y^{max}, {\rm a.e}
\end{equation}
\begin{equation}
\label{lm2}
\forall \delta > 0, \;  {\rm Prob} (Y \in [y^{max} -\delta, y^{max}]) > 0.
\end{equation}
Then we have:
\begin{equation}
\label{lm3}
\lim_{\theta\rightarrow \infty} \frac{\E Y e^{\theta Y}}{\E e^{\theta Y}}= y^{max}.
\end{equation}
and for any continuous real function $G$, we also have:
\begin{equation}
\label{lm4}
\lim_{\theta\rightarrow \infty} \frac{\E G(Y) e^{\theta Y}}{\E e^{\theta Y}}= G(y^{max}).
\end{equation}
\end{lm}
{\it Proof.}
This lemma can be proved by a $(\delta, \varepsilon)$ argument. If $Y$ has a density it can be proved by an integration by parts.$\Box$
\\

\begin{lm}
\label{lemma-PLA2}
Under the assumptions of Proposition \ref{PropPLA1}, the random variable $Y_{n}$ given by (\ref{PLA1}) has an essential maximum $y^{max}=-y^{min}$.
\end{lm}

The remainder of this section is devoted to the proof of lemma \ref{lemma-PLA2}. To this end, we have to prove the conditions of lemma \ref{lemma-MaxRV} are satisfied. The first condition:
\begin{equation}
\label{PLA3}
Y_{n} \leq y^{max}, {\rm a.e}
\end{equation}
is clearly fulfilled. We make a focus now on the second condition
\begin{equation}
\label{PLA4}
\forall \delta_{1} > 0, {\rm Prob} (Y_{n} \geq y^{max} -\delta_{1}) > 0.
\end{equation}
Since the functional $\chi \mapsto B_{n}(\chi)$ is continuous w.r.t the supremum norm $\|\chi\|_{\infty}= \sup_{x}|\chi(x)|$, we have : 
\begin{equation}
\label{PLA5}
\exists \delta : \; \{ \|Y_{n} - \chi^{min}\|_{\infty} \leq \delta \} \subset \{  Y_{n} \geq y^{max} -\delta_{1}  \}
\end{equation}
We shall then prove that :
\begin{equation}
\label{PLA6}
\forall \delta > 0, {\rm Prob} (\|Y_{n} - \chi^{min}\|_{\infty} \leq \delta) > 0.
\end{equation}
The proof of (\ref{PLA6}) relies on some results concerning Gaussian random fields. For completeness we recall their full statements. 

The first one is the small deviations lower bound of Gaussian random processes (Talagrand \cite{Talagrand}), see Ledoux \cite{Ledoux} (p. 257), Li and Shao \cite{Li-Shao} (Theorem 3.8). We consider a Gaussian separable process $X_{t}$ indexed by a set $T$. This set is equipped by a metric $\rho$ usually called the Dudley metric defined by $\rho(s,t)= (E(X(t)-X(s))^{2})^{1/2}$ (although it does not necessarily separate the points of $T$). The entropy numbers or covering numbers $N(T,\rho;\varepsilon)$ are defined as the minimal number of open balls of radius $\varepsilon$ w.r.t. the metric $\rho$ that are necessary to cover $T$, that is $N(T,\rho;\varepsilon)$ is the smallest number $n$ such that there exist $t_{1}, ..., t_{n}\in T$ : for all $t\in T$, we have $t\in B(t_{i}, \rho; \varepsilon)$ (ball of radius $\varepsilon$ in $T$) for some $i=1, ..., n$, i.e. $\rho(t, t_{i}) < \varepsilon$. Then :

\begin{thm}
\label{Talagrand}
Let $X_{t}$ be a separable Gaussian process indexed by a set $T$ with $N(T,d, \varepsilon)$ its entropy number. Suppose that there exists a nonnegative function $\zeta$ and constants $c_{1}, c_{2}$ such that $ N(T,d, \varepsilon) \leq \zeta(\varepsilon)$ and $ 0 < c_{1} \zeta(\varepsilon)\leq \zeta(\varepsilon/2)\leq c_{2} \zeta(\varepsilon) < \infty $. Then, there exists a constant $K> 0$ such that for every $\varepsilon > 0$:
\begin{equation}
\label{Talagrand1}
\log {\rm Prob}(\sup_{s,t\in T}\| X_{s}-X_{t}\| \leq \varepsilon) \geq -K \zeta (\varepsilon)
\end{equation}
\end{thm}
The second result we need is a weaker version the Gaussian correlation inequality (Li \cite{Li99}, cf. also Li and Shao \cite{Li-Shao}, Theorem 2.14):

\begin{thm}
\label{Li}
Let $\mu$ be a Gaussian measure on a separable Banach space $E$. Then for any symmetric convex sets $A, B \subset E$, and $\lambda \in ]0,1[$ we have:
\begin{equation}
\label{Li1}
\mu (A\cap B) \geq \mu (\lambda A) \mu ((1-\lambda^{2})^{1/2} B)
\end{equation}
and for any centered jointly Gaussian variables $X, Y$ we have
\begin{equation}
\label{Li2}
{\rm Prob} (X\in A, Y\in  B) \geq {\rm Prob} (X\in A) {\rm Prob} (Y \in (1-\lambda^{2})^{1/2} B)
\end{equation}
\end{thm}
We could use the celebrated Gaussian correlation inequality: for {\it symmetric} convex sets $A, B \in \R^{p}$
\begin{equation}
\label{Royen}
{\rm Prob} (X\in A, Y\in  B) \geq {\rm Prob} (X\in A) {\rm Prob} (Y \in B) 
\end{equation}
(Royen \cite{Royen}), but we will need only the previous version. Note that (\ref{Royen}) is also valid for Gaussian measures or Gaussian random variables in a Banach space. 

These results will be used within the following framework:

{\bf Settings:} For a finite volume $V\subset\R^{d}$, let us consider the Banach space $E_{c}=C(V)\subset {\cal S}'(\R^{d})$ of continuous functions from $V$ to $\R$ equipped with the supremum norm. We recall also that $n$ is fixed and we consider the random field $\psi_{n}(x):= \varphi_{n}(x)/\sqrt{c_{n}}$ which is a regular Gaussian random field that gives rise to the measure $\nu_{n}$ on $E_{c}$ which is its law.

Now, we focus on the proof of (\ref{PLA6}). 

{\it Step 1:} We first show that :
\begin{equation}
\label{PLA8}
{\rm Prob}(\sup_{x,y\in V_{ni}}| \psi_{n}(x)-\psi_{n}(y)| \geq \eta_{1}) > 0
\end{equation}
This claim is a consequence of the theorem \ref{Talagrand} above, provided we prove that its conditions are satisfied. In our case, since $\E \psi_{n}(x)^{2}=1$, we have:
\begin{equation}
\label{PLA9}
\rho (x,y)^{2}=\E (\psi_{n}(x)-\psi_{n}(y))^{2}= 2-2 \E \psi_{n}(x)\psi_{n}(y)=2(1- \frac{c_{n}(x-y)}{c_{n}})
\end{equation}
We have to estimate the minimal number $N(T,\rho;\varepsilon)$ of open balls of radius $\varepsilon$ w.r.t. the metric $\rho$ that are necessary to cover $T= V_{ni}$. Let us first see what would be a ball $B(x_{0}, \varepsilon)$ centered in $x_{0}\in T$ and with radius $\varepsilon$ with respect to the distance $\rho$. In view of (\ref{PLA9}), we have 
\[\rho(x_{0}, x) \leq \varepsilon \Longleftrightarrow  \frac{c_{n}(x-x_{0})}{c_{n}} \geq 1-\frac{\varepsilon}{2} \]
Let us take $x_{0}=0$ for instance. Then, from the estimate (\ref{Not10b}), we deduce that $c_{n}(x)/{c_{n}} \geq 1-\frac{\varepsilon}{2}$ is equivalent to:
\[   \int_{\R^{d}} \frac{e^{-v^{2}}}{|v-nx|^{d-2}}dv / \int_{\R^{d}} \frac{e^{-v^{2}}}{|v|^{d-2}}dv \geq 1-\frac{\varepsilon}{2}
\]
which is fulfilled if :
\[   |\int_{\R^{d}} \frac{e^{-v^{2}}}{|v-X|^{d-2}} -\frac{e^{-v^{2}}}{|v|^{d-2}}| dv=  
               |\int_{\R^{d}} e^{-v^{2}} \frac{|v-X|^{d-2}-|v|^{d-2}}{|v-X|^{d-2}|v|^{d-2}}| dv   \leq \frac{K\varepsilon}{2} \]
For some $K > 0$, where we have set $X=nx$ and used the fact that $ \int e^{-v^{2}}/|v|^{d-2}< \infty$. Now since 
\begin{align*}
 ||v-X|^{d-2}-|v|^{d-2}| & = |(|v-X|-|v|)(|v-X|^{d-3}+|v-X|^{d-4}|v|+ ... +|v|^{d-3})| \\
                         &\leq |X|(|v-X|^{d-3}-|v-X|^{d-4}|v|+ ... +|v|^{d-3})|
\end{align*}
we get in fact that $\rho(0, x) \leq \varepsilon $ if (and only if, as it can be verified) $|X|=|nx| \leq K_{1} \varepsilon $ for some $K_{1}>0$. Therefore the number of balls of radius $\varepsilon$ with respect to the metric $\rho$ verifies : $N(T,\rho;\varepsilon) \leq \zeta(\varepsilon)= K_{2} |T|/\varepsilon$, with $K_{2}>0 $. 
The function $\zeta(\varepsilon)$ satisfies obviously the condition $ 0 < c_{1} \zeta(\varepsilon)\leq \zeta(\varepsilon/2)\leq c_{2} \zeta(\varepsilon) < \infty $ of Theorem \ref{Talagrand}.

{\it Step 2:} To continue the proof of (\ref{PLA6}), we begin with the case where the minimal path $\chi^{min}$ is constant : $\chi^{min}(x)=Const$. This happens for instance for :
\begin{equation}
\label{PLA18}
B_{n}(\chi)= -\alpha_{n} \int_{V_{ni}}:\chi^{4}(\frac{x}{n})):dx
\end{equation}
in which case the minimum is $y^{min}= -6\alpha_{n} |V_{ni}|$ and it is attained when $\chi\equiv \sqrt{3}$ or $\chi\equiv -\sqrt{3}$. This also happens if $\beta_{n}\neq 0$ and $\gamma_{n}= 0$.

For $\delta > 0 $, we have
\begin{align}
\label{PLA19}
{\rm Prob} (\sup_{x\in V}|\psi_{n}(x) - y^{min}| \leq \delta) & \geq {\rm Prob} (\sup_{x\in V}|\psi_{n}(x) - \psi_{n}(0)| \leq \delta/2, \nonumber \\
                             & \sup_{x\in V}|\psi_{n}(0) - y^{min}| \leq \delta/2).
\end{align}
We wish to apply the Gaussian correlation inequality (in Banach space this time) in order to have:
\begin{align}
\label{PLA20}
{\rm Prob} (\sup_{x\in V}|\psi_{n}(x) & - \psi_{n}(0)| \leq \delta/2, \sup_{x\in V}|\psi_{n}(0) - y^{min}| \leq \delta/2) \geq \nonumber \\
 & {\rm Prob} (\sup_{x\in V}|\psi_{n}(x) - \psi_{n}(0)| \leq \delta/2) {\rm Prob} (\sup_{x\in V}|\psi_{n}(0) - y^{min}| \leq \delta/2) 
\end{align}
There is however a problem : the sets $\{ \sup_{x\in V}|\psi_{n}(x) - \psi_{n}(0)| \leq \delta/2\}$ and $\{\sup_{x\in V}|\psi_{n}(0) - y^{min}| \leq \delta/2) \}$ are convex but the former is symmetric and the later is not. The symmetry property is a basic assumption for the Gaussian correlation inequality ; there are only few results relaxing this assumption (e.g., \cite{DGM}, \cite{Hu}, \cite{SW}). Y.Hu proved another correlation inequality involving convex functions (and not convex sets) :
\begin{equation}
\label{Hu}
\int_{\R^{d}} f(x) g(x) \mu(dx) \geq \int_{\R^{d}} f(x) \mu(dx)\int_{\R^{d}} g(x) \mu(dx) 
\end{equation}
where $f, g$ are convex and one of them is even ($f(-x)=f(x)$). For the correlation inequality (\ref{Royen}) to be equivalent to (\ref{Hu}), the functions $f,g$ need to be quasi-concave. The paper \cite{DGM} considers convex sets and layers and the later ends with a statement of a conjecture for (\ref{Royen}) for two convex sets where one of them in not symmetric about the origin.

In fact in our case we do not necessarily need the inequality (\ref{PLA20}), but we only need to prove that :
\begin{equation}
\label{PLA21}
 {\rm Prob} (\sup_{x\in V}|\psi_{n}(x) - \psi_{n}(0)| \leq \delta/2,  \sup_{x\in V}|\psi_{n}(0) - y^{min}| \leq \delta/2) >  0
\end{equation}
To this end we shall use the quasi-invariance of the Gaussian measure $\nu_{n}$. Let us consider the translation : $T : \psi \mapsto \psi + \tilde{y}^{(min)}$ on $E$, where $\tilde{y}^{(min)}$ denotes function on ${\cal S}$ which equals to $y^{min}$ on the volume $V$. Then, for the set :
\begin{equation}
\label{PLA22}
 B_{\delta}=\{\psi \in E : \sup_{x\in V}|\psi_{n}(x) - \psi_{n}(0)| \leq \delta/2,  \sup_{x\in V}|\psi_{n}(0) - y^{min}| \leq \delta/2) \},
\end{equation}
we have :
\begin{align}
\label{PLA23}
T( B_{\delta})  & =\{\psi \in E : \sup_{x\in V}|\psi_{n}(x) - \psi_{n}(0)| \leq \delta/2, \sup_{x\in V}|\psi_{n}(0)| \leq \delta/2) \}  
\end{align}
Since the sets $ \{\psi \in E : \sup_{x\in V}|\psi_{n}(x) - \psi_{n}(0)| \leq \delta/2 \}$ and $\{ \sup_{x\in V}|\psi_{n}(0)| \leq \delta/2) \}$ are this time convex and symmetric, the Gaussian correlation inequality or the weaker version ({\ref{Li}) can be applied and we get :
\begin{align}
\label{PLA24}
\nu_{n}(T( B_{\delta})) & ={\rm Prob}(\sup_{x\in V}|\psi_{n}(x) - \psi_{n}(0)| \leq \delta/2, \sup_{x\in V}|\psi_{n}(0)| \leq \delta/2) \nonumber \\ 
							 & \geq {\rm Prob} ( \sup_{x\in V}|\psi_{n}(x) - \psi_{n}(0)| \leq \delta/2) {\rm Prob} ( \sup_{x\in V}|\psi_{n}(0)| \leq \delta/2) \nonumber \\
							& > 0 
\end{align}
because ${\rm Prob} ( \sup_{x\in V}|\psi_{n}(x) - \psi_{n}(0)| \leq \delta/2) >  0$ by the (\ref{PLA8}) and ${\rm Prob} ( \sup_{x\in V}|\psi_{n}(0)| \leq \delta/2)$ as $\psi_{n}(0 >  0)$ as $\psi_{n}(0)$  is a standard Gaussian random variables. Now the quasi-invariance property of $\nu_{n}$ implies that $\nu_{n}( B_{\delta})  > 0$ as soon as $\nu_{n}(T( B_{\delta}))  > 0$, so that (\ref{PLA21}) is valid. This completes the proof of (\ref{PLA21}) and (\ref{PLA6}) in the case where $\chi^{min}$ is constant.
\\
\noindent
{\it Step 3:} Proof of (\ref{PLA6}) in the case where the minimal path $\chi^{min}$ is not constant. We may suppose that $\chi$ is regular (see Remark \ref{LA-remark} below). In particular, uniformly continuous, so that there are covering sets $V_{i}, i=1, ..., l$ of $V$ such that :
\begin{equation}
\label{PLA25}
\sup_{x, y\in V_{i}} | \chi^{min}(x) - \chi^{min}(y) | \leq \delta / 8
\end{equation}
Let $x_{i} \in V_{i}$ be some fixed points. Any $x\in V$ belongs to some $V_{i}$ and we have: $|\psi_{n}(x) - \chi^{min}(x) | \leq |\psi_{n}(x) - \chi^{min}(x_{i})|+  |\chi^{min}(x_{i}) - \chi^{min}(x) | \leq  |\psi_{n}(x) - \chi^{min}(x_{i})| + \delta / 8 $ ; hence: 
\begin{align}
\label{PLA26}
\{ \sup_{x \in V_{i}} |\psi_{n}(x) - \chi^{min}(x) | & \leq \delta \} \supset \{ \sup_{x \in V_{i}} |\psi_{n}(x) - \chi^{min}(x_{i}) |  \leq 7 \delta/8 \}   \nonumber \\
                   & \supset \{ \sup_{x \in V_{i}}|\psi_{n}(x) - \psi_{n}(x)(x_{i}) | \leq \delta / 4, \nonumber \\
									& |\psi_{n}(x_{i}) - \chi^{min}(x_{i}) |  \leq \delta/4, i=1, ..., l \} \\
									& =: B
\end{align}
Consider now the translation $T: \psi \mapsto \psi + \chi^{min}$ in the space $E$. By the quasi-invariance of $\nu_{n}$ under $T$, we will have $\nu_{n} (B) > 0  $ i.e. (\ref{PLA6})  as soon as we have $\nu_{n} (T(B)) > 0  $. We have:
\begin{align}
\label{PLA27}
T(B) = &\{ \psi \in E : \sup_{x \in V_{i}}|\psi_{n}(x) + \chi^{min}(x) - (\psi_{n}(x)(x_{i})+\chi^{min}(x_{i})) | \leq \delta / 4, \nonumber \\
        & |\psi_{n}(x_{i}) + \chi^{min}(x_{i}) - \chi^{min}(x_{i}) |  \leq \delta/4, i=1, ..., l   \}
\end{align}
Using : 
\begin{align}
\label{PLA27b}
 |\psi_{n}(x) + \chi^{min}(x) - (\psi_{n}(x_{i})+\chi^{min}(x_{i})) | & \leq   |\psi_{n}(x)- \psi_{n}(x)(x_{i}) | + |\chi^{min}(x) - \chi^{min}(x_{i}) | \nonumber \\
            & \leq |\psi_{n}(x)- \psi_{n}(x_{i}) |+ \delta/8,
\end{align}
 we get:
\begin{align}
\label{PLA28}
T(B)  \supset & \{ \psi \in E : \sup_{x \in V_{i}}|\psi_{n}(x) - \psi_{n}(x_{i})| \leq \delta / 8, \nonumber \\
        & |\psi_{n}(x_{i}) |  \leq \delta/4, i=1, ..., l   \} \nonumber \\
			\subset & \bigcap_{i=1}^{n} A_{i}\bigcap \B_{i}
\end{align}
where : $ A_{i}= \{ \psi \in E : \sup_{x \in V_{i}}|\psi_{n}(x) - \psi_{n}(x_{i})| \leq \delta / 8\} $ and $\B_{i}=\{ |\psi_{n}(x_{i}) |  \leq \delta/4\}$ are this time convex and symmetric sets of $E$. The Gaussian correlation inequality or its weaker form can be applied and we get :
\begin{equation}
\label{PLA28}
\nu_{n}(T(B)) \geq \prod_{1}^{l} \nu_{n}(A_{i}) \nu_{n}(\B_{i}) >  0,
\end{equation}
because $\nu_{n}(A_{i}) >  0$ and $\nu_{n}(B_{i}) >  0$ by the argument of step 2.  This proves (\ref{PLA21}) and (\ref{PLA6}) in the case where $\chi^{min}$ is not constant and completes the proof of Proposition \ref{PropPLA1} and Theorem \ref{ThmPLA}. $\Box$
\begin{Rk}
\label{LA-remark}
In the case where $\chi^{min}$ is not continuous or regular, the previous results are also valid. In fact we may take a norm $\|.\|$ which is weaker than $\|.\|_{\infty}$. We only need a norm that makes the map $\chi \mapsto B_{n}(\chi) $ continuous, we may choose for instance : $\|\chi\|=\|\chi\|_{4} + \|\nabla \chi\|_{2}$, and we adjust the space $E$ accordingly ; we should however require that $\chi^{\min}$ is in this new space. 
\end{Rk}

\section{Weak limit of the interacting field}

In this section we consider the $\varphi^{4}_{d}$ quantum field  model and the renormalized action in a finite volume $V$ :
\begin{equation}
\label{IF1}
{\cal A}_{n}= g_{n} \int_{V}: \varphi_{n}^{4}:+ m_{n} \int_{V}: \varphi_{n}^{2}:+a_{n} \int_{V}: (\partial\varphi_{n})^{2}:
\end{equation}
and the sequence of measures :
\begin{equation}
\label{IF2}
\frac{d \mu_{n}}{d \mu_{0}}= \frac{e^{-{\cal A}_{n}}}{Z_{n}}, \; \; Z_{n} = \E e^{-{\cal A}_{n}}
\end{equation}
The problem is whether there exists a sequence of renormalized constants $ g_{n},  m_{n},  a_{n}$ such that $\mu_{n}$ has a limit in some sense as $n\rightarrow \infty$. In \cite{AA1}, it was proved that such limit does not exist in the strong sense (setwise). We consider now the possibility of its existence in a weak sense. To this end we deal with the issue of the existence or not of the limits of the moments: 
\[ S_{n,p}(f) := \E_{n} \varphi(f)^{p}   \; {\rm or} \; \; S_{n,p}(f^{(1)}, ...., f^{(p)}) := \E_{n} \varphi(f^{(1)}) ... \varphi(f^{(p)}) \]
We shall use the dynamic equation (\ref{MC8b})
\begin{equation}
\label{MC8b1}
\E_{n} \varphi(f) \varphi(\tilde{h})= (h,f) - \E_{n} \varphi(f) \ph1 \A
\end{equation}
and for the action (\ref{IF1}) we have
\begin{equation}
\label{IF3}
\ph1 \A =  \int_{V} dx [ 4 g_{n} :\varphi_{n}^{3}(x): h_{n}(x)+ 2m_{n} \varphi_{n}(x) h_{n}(x)+ 2a_{n}  \sum_{i=1}^{d}\partial_{i}\varphi_{n}(x) \partial_{i} h_{n}(x) ]
\end{equation}
To simplify the last term we use the integration by part formula in $\R^{d}$:
\begin{equation}
\label{IF4}
 \int_{V} u \frac{\partial v}{\partial x_{i}} d^{d}x =  \int_{\partial V} u v n_{i} d\Gamma  -  \int_{V} v \frac{\partial u}{\partial x_{i}} d^{d}x
\end{equation}
with $\partial V$ being the edge surface of the volume $V$ and $n_{i}(x)$ the component $i$ of the normal vector to $\partial V$ at the point $x$.
We apply this to :
\begin{align}
\label{IF5}
\partial_{i}\varphi_{n}(\partial_{i} h_{n})  = & \int_{V} \partial_{i}\varphi_{n} (x)\partial_{i} h_{n}(x) dx    \nonumber \\
                    = &  \int_{\partial V} \varphi_{n} (x)\partial_{i} h_{n}(x)d\Gamma  - \int_{V} \varphi_{n} (x)\partial_{i}^{2} h_{n}(x) dx  \nonumber \\
                    = & - \int_{V} \varphi_{n} (x)\partial_{i}^{2} h_{n}(x) dx  = - \varphi_{n}(\partial_{i}^{2} h_{n}) 
\end{align}
provided we choose $h$ such that $\partial_{i} h_{n}=0$ on the border surface $\partial V$. In the following we shall consider functions $h$ such that $h(x)=Const $ in an arbitrary small neighborhood of $\partial V$, e.g. $h\in {\cal D}(V)$. The formula (\ref{IF3}) becomes:
\begin{equation}
\label{IF6}
\ph1 \A =   4 g_{n} \int_{V} :\varphi_{n}^{3}(x): h_{n}(x) dx+ 2m_{n} \varphi_{n}( h_{n}) - 2a_{n}  \sum_{i=1}^{d}\partial_{i} \varphi_{n}(\partial_{i}^{2} h_{n}) 
\end{equation}
and the equation (\ref{MC8b1}) may be written as:
\begin{align}
\label{IF7}
\E_{n} \varphi(f) \varphi(\tilde{h})= & (h,f)  - \E_{n} [ \varphi(f)  4 g_{n} \int_{V} :\varphi_{n}^{3}(x): h_{n}(x) dx +  2m_{n}\varphi(f) \varphi_{n}( h_{n}) \nonumber \\
                    - &  2 a_{n} \varphi(f) \varphi_{n}(\Delta h_{n}) ] 
\end{align}
We shall also use the second dynamic equation:
\begin{equation}
\label{DynaEq2}
S_{n}^{(2)}(f):= \E_{n}\varphi(f)^{2}= \sf + \E_{n} (\partial_{f} \A_{n})^{2}- \E_{n} \partial_{f}^{2} \A_{n}
\end{equation}
Using the dynamic equation (\ref{MC8})-(\ref{MC8b}) we get the following expression of $\E_{n} (\partial_{f} \A_{n})^{2}$:
\begin{align}
\label{DynaEq2b}
\E_{n} (\partial_{f} \A_{n})^{2}= & \E_{n} \{ 4^{2}g_{n}^{2} [ (\int_{V} \varphi_{n}^{3}f_{n})^{2} - 3^{2} c_{n}^{2} \varphi_{n}(f_{n})^{2}] \nonumber \\
            + & 24 g_{n} c_{n} [ (f,f_{nn}) -  (\varphi(\f)\varphi_{n}(f_{n}) - 2m_{n} \varphi_{n}(f_{n})^{2} + 2a_{n} \varphi_{n}(f_{n}) \varphi_{n}(\partial^{2}f_{n}) ] \\ \nonumber
						- & (2m_{n})^{2} (\varphi_{n}(f_{n}))^{2} - (2a_{n})^{2} (\varphi_{n}( \partial^{2}f_{n}))^{2} \\ \nonumber
						+ & 2(2m_{n}) (2a_{n}) \varphi_{n}(f_{n}) \varphi_{n}( \partial^{2}f_{n}) + 2(2m_{n})( (f,f_{n}) \\ \nonumber
						- & \varphi (\f)\varphi_{n}(f_{n})) + 2(2a_{n})(  \varphi (\f)\varphi_{n}( \partial^{2}f_{n}) - (f,f_{2,n})\}  
\end{align}
The calculus leading to this formula, is a little tedious but straightforward.
\\
\noindent
According to \cite{AA1} and the previous sections we are led to distinguish several situations that are dependent on the dominant terms in the actions ${\cal A}_{n}$.

\subsection{The case where the mass renormalization term is predominant}

We begin by the case where the mass renormalization term is predominant in the action ${\cal A}_{n}$, that is:  $ m_{n} \gg g_{n} c_{n} $ and $ m_{n} \gg  a_{n}$. In fact it will be sufficient to have  $  g_{n} c_{n} \leq m_{n} (1-\varepsilon)/6$ for large $n$ and for some given $\varepsilon > 0 $ arbitrary small, in other words : $\limsup 6g_{n} c_{n}/m_{n} < 1 $. 

\begin{prop}
\label{PropIF1}
In the case 2 where : $  g_{n} c_{n} \leq m_{n} (1-\varepsilon)/6$ for large $n$ and for some given $\varepsilon >$ arbitrary small and $ m_{n} \gg  a_{n}$, we have, for $p\geq 1$ :
\begin{equation}
\label{IF8}
\lim_{n\rightarrow \infty} \E_{n} \varphi(f)^{p} = 0.
\end{equation}
The same result holds for the limit of the $\E_{n} \varphi(f^{(1)}) ... \varphi(f^{(p)}) $, which means that the interacting $\varphi^{4}$ field in this case is trivial in the sense that $\mu_{\infty}=\delta_{0}$
\end{prop}
{\it Proof.} 
We shall prove (\ref{IF8}) for $p=2$, the proof is similar for the general case and for $\E_{n} \varphi(f^{(1)}) ... \varphi(f^{(p)}) $. Let $f\in {\cal S}$; we suppose first that $f\geq 0 $. By applying (\ref{IF7}) to $h=f$ we get:
\begin{align}
\label{IF9}
\E_{n} \varphi(f) \varphi(\tilde{f})= & (f,f)  - \E_{n} [4 g_{n} \varphi(f) \int_{V} \varphi_{n}^{3}(x) f_{n}(x) dx + 12 g_{n} c_{n} \varphi(f) \varphi_{n}( f_{n})  \nonumber \\
                    & -  2m_{n}\varphi(f) \varphi_{n}( f_{n}) +  2 a_{n} \varphi(f) \varphi_{n}(\Delta f_{n}) ] 
\end{align}
By Theorem \ref{Distributions}, we have for the last two terms (if the weak limit of the field is to exist) :
\begin{equation}
\label{IF10}
\E_{n} \varphi(f)  \varphi_{n}(f_{n}) = \E_{n} \varphi(f)^{2} (1+o(n)) \; \; 
          {\rm and} \; \;  \E_{n} \varphi(f)  \varphi_{n}(\Delta f_{n}) = \E_{n} \varphi(f)  \varphi(\Delta f) (1+o(n))
\end{equation}
As for the term containing $\varphi_{n}^{3}(x)$, we note by the first Griffiths inequality for the cutoff field measure $\mu_{n}$ (Proposition \ref{CorrelationIneq}):
\begin{equation}
\label{IF11}
\E_{n} \int_{V} \varphi(f) \varphi_{n}^{3}(x) f_{n}(x) dx = \int_{V} \E_{n} \varphi (\delta_{x}^{n})^{2}(x) \varphi(f) \varphi(\delta_{x}^{n}) f_{n}(x) dx \geq 0
\end{equation}
 which are valid here with the choice of $f\geq 0$, and noting that we always have $\delta_{x}^{n} \geq 0$. (\ref{IF10})-(\ref{IF11}) imply that:
\begin{align}
\label{IF12}
\E_{n} \varphi(f) \varphi(C^{-1}f) &  \leq (f,f) -(2 m_{n} +  12 g_{n} c_{n}) \E_{n} \varphi(f)^{2} + 2 a_{n}  \E_{n} \varphi(f) \varphi(\Delta f)  \nonumber \\
									& + [( m_{n} + g_{n} c_{n}) \varphi(f)^{2} + a_{n}  \E_{n} \varphi(f) \varphi(\Delta f)] o(n)  \nonumber \\
									& \leq (f,f) -(2 m_{n} +  12 g_{n} c_{n}) \E_{n} \varphi(f)^{2} + 2 a_{n}  (\E_{n} \varphi(f)^{2})^{1/2} ( \E_{n} \varphi(\Delta f)^{2})^{1/2}  \nonumber \\
									& + [( m_{n} + g_{n} c_{n})  \E_{n} \varphi(f)^{2} + a_{n}  \E_{n} \varphi(f) \varphi(\Delta f)] o(n)  
\end{align}
from which we get:
\begin{align}
\label{IF13}
\E_{n} \varphi(f)^{2} & \leq  \frac{1}{|2 m_{n} -  12 g_{n} c_{n}|}\{ (f,f) - \E_{n} \varphi(f) \varphi(\f) 
                                          + 2 a_{n}  \E_{n} (\varphi(f)^{2})^{1/2} (\E_{n}\varphi(\Delta f))^{1/2} \nonumber \\
										& + o(n) [( m_{n} + g_{n} c_{n})  \E_{n} \varphi(f)^{2} + a_{n}  \E_{n} (\varphi(f)^{2})^{1/2} (\E_{n}\varphi(\Delta f))^{1/2}] \}
\end{align}
and we see that if a weak limit of the interacting field exists then $\E_{n} (\varphi(h)^{2})$ should be bounded for each $h$ and in the condition of the proposition we have:
\begin{equation}
\label{IF14}
 |2 m_{n} -  12 g_{n} c_{n}|= 2 m_{n} |1- 6  g_{n} c_{n}/ m_{n}| \rightarrow \infty  \; \; {\rm and} \; |2 m_{n} -  12 g_{n} c_{n}| \gg  a_{n},
\end{equation}
(\ref{IF13}) implies that $\E_{n} \varphi(f)^{2} \rightarrow 0$. For a function $f$ which is not positive, we note that we can find $h\in {\cal S}, h\geq 0$ such that $f+h\geq 0$. 
By the previous result we will have $\E_{n} \varphi(f+h)^{2}  \rightarrow 0$ from which we deduce that $\E_{n} \varphi(f)^{2} \rightarrow 0$. Indeed, it suffices to write:
\[ \E_{n} \varphi(f+h)^{2}= \E_{n} \varphi(f)^{2} + \E_{n} \varphi(h)^{2} + 2 \E_{n} \varphi(f)\varphi(h) \]
 and to use the Schwarz inequality and the fact that we already have $\E_{n} \varphi(h)^{2} \rightarrow 0$ for $h\geq 0$.   $\Box$ 
\begin{Rk}
We note that the previous result is valid for all dimensions $d\geq 1$.
\end{Rk}
\subsection{The case where the dominant is the wave renormalization constant : $ a_{n} n^{2}\gg g_{n}c_{n}, m_{n}$ }
\subsubsection {Preliminary results}

In this section we show that the use of the dynamic equations enables to get some first results on the triviality of the $\varphi^{4}$ field in the case where the  wave renormalization constant is predominant. Let us write the first equation $\E_{n} \varphi(f) \varphi(\tilde{h})= (h,f) - \E_{n} \varphi(f) \partial_{h} \A_{n}$ as:
\begin{align}
\label{wave1bis}
2a_{n} \E_{n} \varphi (f) \varphi_{n}(\Delta h_{n}) & = (h,f) - \E_{n}\{ 4g_{n} \varphi (f) \int_{V}: \varphi_{n}^{3}(x): h_{n}(x) \\ \nonumber
                                        &+ 2m_{n} \varphi (f) \varphi_{n}(h_{n}) -\varphi(f) \varphi(\tilde{h} \} 
\end{align}
where we suppose that that $\nabla h_{n}=0$ on $\partial V$. We expect that if $a_{n}$ dominates the terms $m_{n}, g_{n}c_{n}$ we will have $\E_{n} \varphi (f) \varphi_{n}(\Delta h_{n}) \rightarrow 0$ which also means that $\E_{n} \varphi (f) \varphi(\Delta h ) \rightarrow 0$ by Theorem \ref{Distributions}. Also, an inspection of the second dynamic equation (\ref{DynaEq2})-(\ref{DynaEq2b}) makes plausible the fact that $\E_{n} \varphi(\Delta h )^{2} \rightarrow 0 $ when $a_{n}$ is predominant. This is indeed the case, and we have the following:
\begin{prop}
\label{PropWave1}
In the case $ a_{n} n^{2}\gg g_{n}c_{n}, m_{n}$, and under the additional conditions $a_{n} \geq K n^{2} c_{n}g_{n}, K>0$ and $a_{n}\gg m_{n}$, we have for all $f, h \in {\cal S}$, with $f, g \geq 0$ and $\nabla h (x) = 0$ for $x \in \partial V$: 
\begin{equation}
\label{wave1ter}
\E_{n} \varphi (f) \varphi(\Delta h ) \rightarrow 0 \; {\rm as} \; n \rightarrow \infty 
\end{equation}
and
\begin{equation}
\label{wave1quatro}
\E_{n} \varphi(\Delta h )^{2} \rightarrow 0 \; {\rm as} \; n \rightarrow \infty 
\end{equation}
The last limit implies that, if $\mu_{n}$ has a weak limit $\mu_{\infty}$, then the resulting field is trivial in the sense that $\Delta \varphi(f)=0, \; \mu_{\infty}$-a.e.
\end{prop}
Before proving this proposition, we state and prove the following lemma which will be used in this case and the case where the coupling constant is predominant.
\begin{lm}
\label{lemmaWave1}
Let $\eta_{n}: V \subset \R^{d} \rightarrow \R^{+}$ be a sequence of integrable functions such that:
 \[ \forall \varepsilon  > 0, \; \exists n_{\varepsilon} \geq 1 \; \forall n\geq n_{\varepsilon} : \; \; \;  \int_{V} \eta_{n} \leq \varepsilon |V|  \]
Then, there exist a volume $\tilde{V}\subset V$ and a subsequence $n_{k}$ such that : for every $\varepsilon  > 0$, there is $l\geq 1$ such that $\eta_{n_{k}}(x) \leq \varepsilon$ for all $x \in \tilde{V}_{l}$ and $n_{k}\geq n_{l}$. Furthermore $V_{k}\subset V_{l}$ if $k\leq l$ and we have $V_{l}\uparrow V$ as $l\rightarrow \infty$, 
with $|V_{l}| \geq (1-2\sqrt{\varepsilon})|V|$
\end{lm}
{\it Proof.} First, let $\eta : V \subset \R^{d} \rightarrow \R^{+}$ be a function such that $ \int_{V} \eta \leq \varepsilon |V| $ and for $\delta > 0 $ let us set :
$V_{\delta}=\{ x \in V : \eta(x) \leq \delta \}$ and $V_{\delta}^{c}=\{ x \in V : \eta(x) > \delta \}$. Then for $\delta = \sqrt{\varepsilon}$ we have: $ |V_{\sqrt{\varepsilon}}^{c}| \leq \sqrt{\varepsilon} |V|$. 
\\
\noindent
Now suppose we have $\int_{V} \eta_{n} \leq \varepsilon |V|$ for $n\geq n_{\varepsilon}$ and the sequence of functions $\eta_{n}$. Let first us take $\varepsilon_{k}= 1/2^{2k}$; then, using the previous notation and argument, we will have $|V^{c}_{\varepsilon_{k}}| \leq |V|/2^{k}$ for $n\geq n_{k}:=n_{\varepsilon_{k}}$, which means that $\eta_{n_{k}} \leq \varepsilon_{k}$ outside $V^{c}_{\varepsilon_{k}}$. Let $l \geq 1 $, then for $k \geq l$ we have $\eta_{n_{k}} \leq \varepsilon_{k} \leq \varepsilon_{l}$ outside $\tilde{V}_{l}^{c}:= \bigcup_{k\geq l} V^{c}_{\varepsilon_{k}}$, and $| \tilde{V}_{l}^{c}|\leq \sum_{k\geq l} | V^{c}_{\varepsilon_{k}} | \leq \sum_{k\geq l} 1/2^{k}=1/2^{l-1} |V|$. The volume $ \tilde{V}_{l}^{c}$ can therefore be made arbitrarily small if $l$ is sufficiently large and the $n_{l}$ provides the desired subsequence of the lemma. 
Finally for an arbitrary $\varepsilon > 0 $ we consider the smallest $k\geq 1$ such that $\varepsilon_{k}= 1/2^{2k} \varepsilon$ and we get the same conclusion as for the specific $\varepsilon_{k}$ considered before. $\Box$
\\
\\
\noindent
{\it Proof of Proposition \ref{PropWave1}.} 
We note that (\ref{wave1ter}) and (\ref{wave1ter}) are always valid when $g_{n}=m_{n}=0$ or when $g_{n}=0$ and $a_{n} \gg m_{n}$ a fact that is not obvious from the expression of $\mu_{n}$ but is a straightforward consequence of (\ref{wave1bis}) and (\ref{DynaEq2})-(\ref{DynaEq2b}). The interaction term makes things a little more complicated.

In the following we shall also suppose that $f, h \geq 0$. Then, for the present $\varphi^{4}_{d}$ model we have by the first Griffiths inequality $ \E_{n} \varphi (f) \varphi_{n}^{3}(x)  h_{n}(x) \geq 0 $ and (\ref{wave1bis}) gives:
\begin{equation}
\label{wave2}
\E_{n} \varphi(f) \varphi(\tilde{h}) \leq (h,f) + \E_{n} \{  12 g_{n}c_{n} \varphi (f) \varphi_{n}(h_{n}) - 2m_{n} \varphi (f) \varphi_{n}(h_{n}) +2a_{n} \varphi (f) \varphi_{n}(\Delta h_{n})  
						         \} 
\end{equation}
or:
\begin{equation}
\label{wave3}
2a_{n}  \E_{n} \varphi (f) \varphi_{n}(\Delta h_{n})  \geq (h,f) + \E_{n} [(2m_{n}- 12 g_{n}c_{n}) \varphi (f) \varphi_{n}(h_{n}) - \varphi(f) \varphi(\tilde{h})]
\end{equation}
We shall now seek an upper bound for $\E_{n} \varphi (f) \varphi_{n}(\Delta h_{n})$. 
In this case the interaction measure is :
\begin{equation}
\label{wave3bis}
d\mu_{n} = \exp -a_{n}n^{2} {\cal A}_{n}^{(1)} (\psi_{n}) d\mu_{0}
\end{equation}
where the regularized action has been written as:
\begin{align}
\label{wave4}
{\cal A}_{n}(\varphi)= & a_{n}c_{n}n^{2} \int_{V}: [ \frac{ g_{n} c_{n}}{a_{n}n^{2}} :\psi_{n}^{4}:+ \frac{m_{n}}{a_{n}n^{2}} : \psi_{n}^{2}: 
                                                + \frac{1}{d_{n}} \int_{V}: (\partial\psi_{n})^{2}:  \\ \nonumber
						= & a_{n}c_{n}n^{2} {\cal A}_{n}^{(1)} (\psi_{n}) 
\end{align}
The principle of the least action implies that for any $\varepsilon >0 $, there is some $n_{1} \geq 1$ such that:
\begin{equation}
\label{wave5}
\E_{n} \varphi (f) \varphi(h)  \sim \E_{n} \varphi (f) \varphi(h) 1_{\Sigma_{\varepsilon, n}}, \; \; \Sigma_{\varepsilon,n}=\{ -{\cal A}_{n}^{(1)} (\varphi) \geq (1 - \varepsilon) |V|\}
\end{equation}
which means that we can work on the events $\Sigma_{n, \varepsilon}$ as far as we are calculating expectations under $\E_{n}$ for large $n$. We may then assume that:
\begin{equation}
\label{wave6}
-\int_{V}\frac{ g_{n} c_{n}}{a_{n}n^{2}} :\psi_{n}^{4}: - \int_{V} \frac{m_{n}}{a_{n}n^{2}} : \psi_{n}^{2}: 
                                                - \int_{V} ( \frac{(\partial\psi_{n}(x))^{2}}{n^{2}}-1) \geq (1- \varepsilon) |V|
\end{equation}
{\it Remark.} Without using the principle of the least action, we note that for large $n$, the measure $\mu_{n}$ is supported in the sets 
\[ \Sigma_{n, \varepsilon}^{'}= \{ {\cal A}_{n}^{(1)} (\psi_{n}) \geq 0 \} \]
and we can continue with the same arguments, except that $\varepsilon$ will be just a positive number and not an arbitrarily small number, so that we will get a slightly weaker result.
(\ref{wave6}) may be written as:
\begin{equation}
\label{wave6b}
\int_{V}\frac{ g_{n} c_{n}}{a_{n}n^{2}} \psi_{n}^{4} - 6 \frac{ g_{n} c_{n}}{a_{n}n^{2}} \int_{V} \psi_{n}^{2} + \frac{m_{n}}{a_{n}n^{2}} \int_{V} (\psi_{n}^{2}-1) 
                                                +\int_{V} ( \frac{(\partial\psi_{n}(x))^{2}}{n^{2}} dx 
																								\leq (\varepsilon  -3 \frac{ g_{n} c_{n}}{a_{n}n^{2}}) V
\end{equation}
Since $g_{n} c_{n}/(a_{n}n^{2}) \rightarrow 0$ and $m_{n}/(a_{n}n^{2}) \rightarrow 0$, and taking into account $ \int_{V} \psi_{n}^{2} \leq V^{1/2} (\int_{V} \psi_{n}^{4})^{1/2}   $ the later inequality (\ref{wave6b}) holds only if :
\begin{equation}
\label{wave7}
\int_{V}\frac{ g_{n} c_{n}}{a_{n}n^{2}} \psi_{n}^{4} \leq \varepsilon  \; \; {\rm or } \; \int_{V} \psi_{n}^{4} \leq O(1)
\end{equation}
and 
\begin{equation}
\label{wave7b}
\int_{V} \frac{(\partial\psi_{n}(x))^{2}}{n^{2}} dx  \leq \varepsilon 
\end{equation}
for $n \geq n_{1} $ sufficiently large. We can use either (\ref{wave7}) or (\ref{wave7b}) to get an estimate of $\E_{n}\psi_{n}^{2}(x)$. We continue the proof by using (\ref{wave7}) (see Remark \ref{Remark-Wave1} below for the usage of (\ref{wave7b})) :
\begin{equation}
\label{wave8}
\int_{V}\psi_{n}^{4} \leq \max (O(1), \varepsilon \frac{a_{n}n^{2}}{ g_{n} c_{n}}) \; \; 
{\rm and } \; \int_{V}\psi_{n}^{2} \leq \max (O(1), n\sqrt{\varepsilon} \sqrt{\frac{a_{n}}{ g_{n} c_{n}}})
\end{equation}
for $n \geq n_{1}$. By Lemma \ref{lemmaWave1}, there is a volume $\tilde{V_{\varepsilon}}\subset V$ such that :
\begin{equation}
\label{wave8b}
\E_{n} \psi_{n}^{2} (x) \leq \max (O(1), n\sqrt{\varepsilon} \sqrt{\frac{a_{n}}{ g_{n} c_{n}}}), \; \; \forall x \in \tilde{V_{\varepsilon}},  n \geq n_{1}
\end{equation}
or 
\begin{equation}
\label{wave9}
\E_{n}  \varphi_{n}^{2} (x) \leq \max (O(c_{n}), n\sqrt{\varepsilon} \sqrt{\frac{a_{n}c_{n}}{ g_{n} }}), \; \; \forall x \in \tilde{V_{\varepsilon}},  n \geq n_{1}
\end{equation}
Now (\ref{wave1bis}) implies that :
\begin{align}
\label{wave10}
2a_{n} |\E_{n} \varphi (f) \varphi_{n}(\Delta h_{n})| & \leq |\E_{n} 4g_{n} \int_{V} \varphi (f) \varphi_{n}^{3}(x) h_{n}(x)|  + |\E_{n} \varphi(f) \varphi(\tilde{h})| + |(h,f)|  \\ \nonumber
                        & + 12 g_{n}c_{n} |\varphi (f) \varphi_{n}(h_{n})| +  2m_{n} |\E_{n} \varphi (f) \varphi_{n}(h_{n})|  
\end{align}
We have to control the 1st term of the r.h.s. Suppose that $h\geq 0, f\geq 0$, then:
\begin{align}
\label{wave11}
|\E_{n} 4g_{n} \int_{V} \varphi (f) \varphi_{n}^{3}(x) h_{n}(x) dx |  &= 4g_{n} \int_{V} \E_{n} ( \varphi_{n}^{2}(x) \varphi (f) \varphi_{n}(x) h_{n}(x)) dx  \\ \nonumber
                            & \leq 4g_{n} \int_{V} 3 \E_{n} \varphi_{n}^{2}(x) \E_{n} \varphi (f) \varphi_{n}(x) h_{n}(x)) dx  \nonumber                   
\end{align}
where we used the 1st Griffiths inequality and the Gaussian inequalities. If we suppose that $f,h$ have their supports in $\tilde{V_{\varepsilon}}$, we get by (\ref{wave9}):
\begin{align}
\label{wave12}
|\E_{n} 4g_{n} \int_{V} \varphi (f) \varphi_{n}^{3}(x) h_{n}(x) dx |  & \leq 12 g_{n} n\sqrt{\varepsilon} \sqrt{\frac{a_{n}c_{n}}{ g_{n} }} \int_{V} \E_{n} (\varphi (f) \varphi_{n}(x) h_{n}(x)) dx  \\ \nonumber
               & \leq 12 n\sqrt{\varepsilon} \sqrt{a_{n}c_{n}g_{n} } \E_{n} \varphi (f) \varphi_{n}(h_{n}) \\
							& \leq 12 n\sqrt{\varepsilon} \sqrt{a_{n}c_{n}g_{n} } (\E_{n} \varphi (f) \varphi(h) + o(n))
\end{align}
So that (\ref{wave10}) becomes:
\begin{align}
\label{wave13}
 |\E_{n} \varphi (f)\varphi_{n}(\Delta h_{n})| & \leq  6 n\sqrt{\varepsilon} \sqrt{\frac{c_{n}g_{n}}{a_{n}} } (\E_{n} \varphi (f) \varphi(h) + o(n)) \\ \nonumber
                      &  + \frac{1}{a_{n}} [(|\E_{n} \varphi(f) \varphi(\tilde{h})| + |(h,f)|)  
                       + (6 g_{n}c_{n}  +  m_{n}) |\E_{n} \varphi (f) \varphi_{n}(h_{n})|  ]
\end{align}
By Theorem \ref{Distributions}, we have the same inequality for $ |\E_{n} \varphi (f)\varphi(\Delta h)|$, from which we deduce that if $a_{n} \geq K n^{2} c_{n}g_{n} $ for some $K>0$, then we will have $|\E_{n} \varphi (f)\varphi(\Delta h)| \rightarrow 0$. We have of course the same result if $h\leq 0$. 

We turn now to the proof of (\ref{wave10}). We use to the this end the second dynamic equation (\ref{DynaEq2}) and the extended expression (\ref{DynaEq2b}) of $\E_{n} (\partial_{f} {\cal A}_{n})^{2}$, which give:
\begin{align}
\label{wave14}
(2a_{n})^{2}  \E_{n} (\varphi_{n}( \Delta f_{n}))^{2} & = \E_{n}\varphi(f)^{2} + \sf 
                                    + \E_{n} \{ 4^{2}g_{n}^{2} [ (\int_{V} \varphi_{n}^{3}f_{n})^{2} - 3^{2} c_{n}^{2} \varphi_{n}(f_{n})^{2}]  \\ \nonumber
            + & 24 g_{n} c_{n} [ (f,f_{nn}) -  (\varphi(\f)\varphi_{n}(f_{n}) - 2m_{n} \varphi_{n}(f_{n})^{2} + 2a_{n} \varphi_{n}(f_{n}) \varphi_{n}(\partial^{2}f_{n}) ] \\ \nonumber
						- & (2m_{n})^{2} (\varphi_{n}(f_{n}))^{2} + 2(2m_{n}) (2a_{n}) \varphi_{n}(f_{n}) \varphi_{n}( \partial^{2}f_{n})  \\ \nonumber
						+ &  2(2m_{n})( (f,f_{n}) - \varphi (\f)\varphi_{n}(f_{n})) + 2(2a_{n})(  \varphi (\f)\varphi_{n}( \partial^{2}f_{n}) - (f,f_{2,n})\} \\ \nonumber
						- &  12 g_{n} \E_{n}\int_{V} :\varphi_{n}^{2}:f_{n}^{2}
\end{align}
We want to show that $\E_{n} (\varphi_{n}( \Delta f_{n}))^{2} \rightarrow 0$. As we assume that the limiting field exists, i.e. $ \E_{n} \varphi (f) \varphi (h)$ exists for all $f,h$, all the terms of the r.h.s of (\ref{wave14}) are clearly negligible w.r.t. $a_{n}^{2} $ by the assumption of the Proposition ($a_{n} \geq K n^{2} c_{n}g_{n}$ and $a_{n}\gg m_{n}$), except the first term 
\begin{equation}
\label{wave14b}
t_{1,n} := g_{n}^{2}\E_{n} (\int_{V} \varphi_{n}^{3}f_{n})^{2}= 4^{2}g_{n}^{2}\E_{n} \int_{V}\int_{V} \E_{n} (\varphi_{n}^{3}(x) \varphi_{n}^{3}(y))f_{n}(x) f_{n}(y) dx dy
\end{equation}
 which we will consider now. By the Gaussian and the Schwarz inequalities we have:
\begin{align}
\label{wave14c}
\E_{n} (\varphi_{n}^{3}(x) \varphi_{n}^{3}(y)) &\leq 9 \E_{n} \varphi_{n}^{2}(x) \E_{n}  \varphi_{n}^{2}(y) \E_{n}  (\varphi_{n}(x) \varphi_{n}(y))
                                                   + \E_{n}  (\varphi_{n}(x) \varphi_{n}(y))^{3} \\ \nonumber
                                         & \leq 10 \E_{n} \varphi_{n}^{2}(x) \E_{n}  \varphi_{n}^{2}(y) \E_{n}  (\varphi_{n}(x) \varphi_{n}(y))
\end{align}
(recall that all the previous expectations are positive by the 1st Griffiths inequality). On account of (\ref{wave9}) we can control $t_{1,n}$ as:
\begin{align}
\label{wave14d}
t_{1,n} & \leq 10 g_{n}^{2} (\max_{x} \E_{n}  \varphi_{n}^{2}(x))^{2} \int_{V}\int_{V} \E_{n}  (\varphi_{n}(x) \varphi_{n}(y)) f_{n}(x) f_{n}(y) dx dy \\ \nonumber
        & \leq 10 g_{n}^{2}  \varepsilon n^{2} \frac{a_{n}c_{n}}{ g_{n} }  \E_{n}  (\varphi_{n}(f))^{2}
\end{align}
 which implies that $t_{1,n}/a_{n}^{2} \leq 10 \varepsilon n^{2} g_{n} c_{n}/ a^{n} $. So that $t_{1,n}/a_{n}^{2} \rightarrow 0 $ and by (\ref{wave14}) 
\begin{equation}
\label{wave14e}
|\E_{n} \varphi_{n}(\Delta f)^{2}| \rightarrow 0
\end{equation}
 under the assumptions of Proposition \ref{PropWave1}.
This implies that $\E_{\infty} \varphi(\Delta f)^{2}= 0$ and therefore $\Delta \varphi=0$, $\mu_{\infty}$-a.e., if we assume that $\mu_{n}\rightarrow \mu_{\infty}$ weakly. We have indeed : $\Delta \varphi_{n} (f_{n})=\varphi_{n} (\Delta f_{n})=\varphi(\Delta f_{n})_{n} $. Noting that if $h_{n} \rightarrow  h$ then $h_{nn} \rightarrow  h$ for any $h\in {\cal S}$, we also have : $\Delta f_{n} \rightarrow \Delta f$ and $(\Delta f_{n})_{n} \rightarrow \Delta f$. This and (\ref{wave14e}) show that $\E_{\infty} \varphi(\Delta f)^{2}= 0$ by Theorem \ref{Distributions}. $\Box$
\begin{Rk}
\label{Remark-Wave1}
We can also use (\ref{wave7b}) or the fact that $(\partial \psi_{n})(x)^{2} \leq K n$ for some $K\geq 0$ (or even $(\partial \psi_{n})(x)^{2} = o(n)$, by using the principle of the least action) on the sets of paths $\Sigma_{n,\varepsilon}\subset {\cal P}$. By using lemma \ref{lemmaWave1} and an argument presented in the Corollary \ref{CoroWave1}  below (which allows to pass from an estimate of $\E_{n} (\partial \varphi_{n})(x)^{2}$ to an estimate of $\E_{n} \varphi_{n}(x)^{2}$), this will lead to the fact that $ \E_{n} \varphi_{n}(x)^{2} = O ( c_{n}n^{2}$ or even $ \E_{n} \varphi_{n})(x)^{2} = o (c_{n}n^{2})$ for all $x\in V$ which is equivalent to (\ref{wave9}) and the remaining of the proof is simialr to the previous one.
\end{Rk}
\begin{Rk}
One expects that $|\E_{n} \varphi (f)\varphi( h)| \rightarrow 0$ for all $f, h\geq 0 $. However, this can not be deduced from Proposition \ref{PropWave1}: if we want to use this proposition, we have to look for a function $u$ such that $\Delta u = h$, but we have other constraints : the system 
\begin{equation}
\label{wave14e}
\begin{cases}
\Delta u = h\\
\nabla u (x) = 0 \; {\rm for} \; x \in \partial V \\
u \geq 0  \; {\rm or} \; u\leq 0 \; 
\end{cases}
\end{equation}
and by the Gauss formula we have $\int_{V} \Delta u = \int_{\partial V} \nabla u $ ; this gives the compatibility condition : $\int_{V} h = 0$. If $h\geq 0$, the previous system has no  solution unless $h=0$. Proving directly that $|\E_{n} \varphi (f)\varphi( h)| \rightarrow 0$ is possible (see below) but seems to require other investigations. 
\end{Rk}
\subsubsection{Triviality of the field in the case where the wave renormalization term is predominant}
The purpose of this section is to show the field $\varphi^{4}_{d}, d\geq 4$ is trivial when the wave renormalization term is predominant. We begin by extending the validity conditions of Proposition \ref{PropWave1}:
\begin{prop}
\label{PropWave1b}
In the case $ a_{n} n^{2}\gg g_{n}c_{n}, m_{n}$, and under the additional conditions $a_{n} \gg c_{n}g_{n}$ and $a_{n}\gg m_{n}$, then the conclusions of Proposition \ref{PropWave1} are valid, in particular if $\mu_{n}$ has a weak limit $\mu_{\infty}$, then the resulting field is trivial in the sense that $\Delta \varphi (f)=0, \; \mu_{\infty}$-a.e.
\end{prop}
We will show the triviality result is valid in the most general case ($ a_{n} n^{2}\gg g_{n}c_{n}, m_{n}$), more precisely, we have the following:
\begin{prop}
\label{PropWave1c}
When the wave renormalization term is predominant in the renormalized Lagrangian, i.e. $ a_{n} n^{2}\gg g_{n}c_{n}, m_{n}$, if one supposes that the sequence of the field measures $\mu_{n}$ has a weak limit $\mu_{\infty}$, then this limit is trivial in the sense $\mu_{\infty}=\delta_{0}$, in other words for all $f\in {\cal S}$ we have: $\varphi(f)=0$, $\mu_{\infty}$- almost everywhere.
\end{prop}
The starting point is the expression of the interaction measure (\ref{wave3bis}) : $d\mu_{n} = \exp (-a_{n}n^{2} {\cal A}_{n}^{(1)}) (\psi_{n}) d\mu_{0}$. From this expression, it is clear that the more $|\nabla \psi_{n}|$ is small, the more the exponent factor of the measure is big. The idea is to see whether the measure $\mu_{n}$ is supported in the sets $|\nabla \psi_{n}|$ is very small, i.e. when the exponent in (\ref{wave4}) is the largest possible. 
\\
Proposition \ref{PropWave1b} will be a consequence of the following propositions :
\begin{prop}
\label{PropWave3}
For $n\geq 1$ fixed and all $\eta > 0 $, there exists $K > 0$ : 
\begin{equation}
\label{wave23}
{\rm Pr}( \max_{x\in V}|\nabla \psi_{n}(x) - \nabla \psi_{n}(0) | \leq \eta ) \geq e^{\frac{-K |V| n^{3}}{\eta}}
\end{equation}
\end{prop}
{\it Proof.} (\ref{wave23}) is a consequence of the Talagrand small deviations lower bound of Gaussian random processes (Theorem \ref{Talagrand}). We can verify the conditions of this theorem as follows. As we consider here the process $\p \psi_{n}(x):= \nabla \psi_{n}(x)$, the distance $\rho $ is :
\begin{equation}
\label{wave23a}
\rho(x,y) = \E (\p \psi_{n}(x) - \p \psi_{n}(y))^{2}= c_{n}^{-1}\E (\p \varphi_{n}(x)^{2} + \p \varphi_{n}(y)^{2}- 2 \varphi_{n}(x)\varphi_{n}(y))
\end{equation}
We note that the expectation is taken w.r.t. $\mu_{0}$. The estimates of the free covariance (\ref{Not10b}) can be written as:
\begin{equation}
\label{wave23b}
\E \varphi_{n}(x) \varphi_{n}(y) \sim K n^{d-2} \int_{\R^{d}} \frac{e^{-(u+n(x-y))^{2}}}{u^{d-2}}du
\end{equation}
and the covariance of the gradient of the field will be estimated by:
\begin{equation}
\label{wave23c}
\E \p\varphi_{n}(x) \p \varphi_{n}(y) \sim 4K n^{d} \int_{\R^{d}} (u+n(x-y))^{2} \frac{e^{-(u+n(x-y))^{2}}}{u^{d-2}}du
\end{equation}
and, using the invariance by translation of $\mu_{0}$, the distance $\rho$ will be given by:
\begin{align}
\label{wave23d}
\rho(x,y)  & =  2 c_{n}^{-1}\E (\p \varphi_{n}(0)^{2} - 2 \varphi_{n}(0)\varphi_{n}(x-y))  \\ \nonumber
           & \sim \frac{K_{1}}{n^{d-2}} n^{d} \int_{\R^{d}}  \frac{u^{2} e^{-u^{2} - (u+n(x-y))^{2} e^{-(u+n(x-y))^{2}}}  }{|u|^{d-2}}du \\ \nonumber
					 & \sim K_{1}n^{2} \int_{\R^{d}} [ \frac{u^{2} e^{-u^{2}}}{|u|^{d-2}} -  \frac{u^{2} e^{-u^{2}}}{|u -n(x-y)|^{d-2}}]  du \\ \nonumber
					 & \sim K_{1}n^{2} \int_{\R^{d}} [ \frac{u^{2} e^{-u^{2}} ||u|^{d-2}  - |u-X|^{d-2}|}{|u|^{d-2}|u-X|^{d-2}}  du ] \\ \nonumber
					 & \leq K_{1}n^{2} |X| \int_{\R^{d}} [ \frac{u^{2} e^{-u^{2}} ||u|^{d-3} + ... |}{|u|^{d-2}|u-X|^{d-2}}  du ] 
\end{align}
where $X=n(x-y)$. The last integral being convergent, we get:
\begin{equation}
\label{wave24}
\rho(x,y) \leq K_{2} n^{3}|x-y|
\end{equation}
Hence, in order to have $\rho(x,y) \leq \varepsilon$ we should have $|x-y| \leq \varepsilon/(Kn^{3}) $ and from this we deduce that the number $N(\varepsilon)$ of balls of radius $\leq \varepsilon$ that are necessary to cover $V$ verifies $N(\varepsilon) \sim V/ (\varepsilon/n^{3}) = V n^{3}/ \varepsilon$. So that with $\zeta (\varepsilon) = N(\varepsilon)$ the conditions of Theorem \ref{Talagrand} are satisfied and there is a $K \geq 0$ such that (\ref{wave23}) is valid. $\Box$
\begin{prop}
\label{PropWave4}
For $n\geq 1$ fixed and all $\eta, \eta^{'} > 0 $, there exists $K > 0$ such that:
\begin{equation}
\label{wave25}
\mu_{n}(\Sigma (\eta, n)) = \mu_{n} \{|\nabla \psi_{n}(x)| \leq \eta, \forall x\in V\} \geq \frac{1}{Z_{n}} e^{a_{n}n^{2} c_{n}(1-\eta |V| /n^{2})} e^{\frac{-K |V| n^{3}}{\eta}}
\end{equation}
\begin{equation}
\label{wave26}
\mu_{n}(\Sigma_{I} (\eta, n)) = \mu_{n} \{\int_{V}|\nabla \psi_{n}(x)| \leq \eta \} \geq \frac{1}{Z_{n}} e^{a_{n}n^{2} c_{n}(1-\eta /n^{2})} e^{\frac{-K |V| n^{3}}{\eta}}
\end{equation}
\begin{equation}
\label{wave27}
\mu_{n}(\Sigma_{I}^{c} (\eta^{'}, n)) = \mu_{n} \{\int_{V}|\nabla \psi_{n}(x)| \geq \eta^{'} \} \leq \frac{1}{Z_{n}} e^{a_{n}n^{2} c_{n}(1-\eta^{'} /n^{2})}  
\end{equation}
and for any (measurable) subset $V_{1}\subset V$, there exists $K > 0$ such that:
\begin{equation}
\label{wave28}
\mu_{n}(\Sigma_{1}^{c} (\eta^{'}, n)) = \mu_{n} \{\int_{V}|\nabla \psi_{n}(x)| \geq \eta^{'}, \forall x\in V_{1} \} 
                                 \leq \frac{1}{Z_{n}} e^{a_{n}n^{2} c_{n}(1-\eta^{'}|V_{1}|  /n^{2})} 
\end{equation}
These inequalities are also valid when we replace $\int_{V}|\nabla \psi_{n}(x)|$ by $\int_{V}|\nabla \psi_{n}(x)|^{2}$.
\end{prop}
{\it Proof.} For the inequality (\ref{wave25}), we have:
\begin{equation}
\label{wave29}
\mu_{n} \{|\nabla \psi_{n}(x)| \leq \eta, \forall x\in V\} \geq \frac{1}{Z_{n}} e^{a_{n}n^{2} c_{n}(1-\eta |V| /n^{2})} 
                    \mu_{0} \{|\nabla \psi_{n}(x)| \leq \eta, \forall x\in V\}
\end{equation}
and 
\begin{equation}
\label{wave30}
\mu_{0} \{\max_{x \in V}|\nabla \psi_{n}(x)| \leq \eta\} \geq \mu_{0} \{ \max_{x \in V}|\nabla \psi_{n}(x) - \nabla \psi_{n}(0)| \leq \eta/2, |\nabla \psi_{n}(0)| \leq \eta/2\}.
\end{equation}
We can then use either the correlation inequalities of Theorem \ref{Li} or the Gaussian correlation inequality (\ref{Royen}): the sets $\{ \max_{x \in V}|\nabla \psi_{n}(x) - \nabla \psi_{n}(0)| \leq \eta/2\}$ and $\{ |\nabla \psi_{n}(0)| \leq \eta/2\}$ are indeed symmetric convex as a set of $E_{c}=C(V)\subset {\cal S}'(\R^{d})$ of continuous functions from $V$ to $\R$ equipped with the supremum norm (the settings are those adopted for the application of Theorems \ref{Talagrand} and \ref{Li} above). Hence we get from (\ref{wave30}) and the Gaussian correlation inequality for instance:
\begin{equation}
\label{wave31}
\mu_{0} \{\max_{x \in V}|\nabla \psi_{n}(x)| \leq \eta\} \geq \mu_{0} \{ \max_{x \in V}|\nabla \psi_{n}(x) - \nabla \psi_{n}(0)| \leq \eta/2\} 
                                                 \mu_{0} \{ |\nabla \psi_{n}(0)| \leq \eta/2\}
\end{equation}
On the other hand since $\nabla \psi_{n}(0)$ is a Gaussian random variable of variance $n^{2}$, we have
\begin{equation}
\label{wave32}
\mu_{0} \{ |\nabla \psi_{n}(0)| \leq \eta/2\} = \int_{-\eta/2}^{\eta/2}\frac{1}{\sqrt{2\pi}n} e^{-\frac{x^{2}}{2n^{2}}} = O(1)
\end{equation}
so that (\ref{wave25}) is a consequence of (\ref{wave31}), (\ref{wave32}) and Proposition \ref{PropWave3}.
The inequality (\ref{wave26}) follows from (\ref{wave31}), Proposition \ref{PropWave3} and the fact that
\begin{equation}
\label{wave33}
\mu_{n} \{\int_{V}|\nabla \psi_{n}(x)| \leq \eta \} \geq \frac{1}{Z_{n}} e^{a_{n}n^{2} c_{n}(1-\eta /n^{2})} \mu_{0} \{\int_{V}|\nabla \psi_{n}(x)| \leq \eta \}
\end{equation}
and $ \{\max_{x\in V}|\nabla \psi_{n}(x)| \leq \eta/|V| \} \subset \{\int_{V}|\nabla \psi_{n}(x)| \leq \eta \} $ and 
The inequalities (\ref{wave27}) and (\ref{wave28}) are straightforward.$\Box$
\begin{prop}
\label{PropWave5}
For any $\eta > 0 $, there exists $ N(\eta) \geq 1$ such that 
\begin{equation}
\label{wave34}
\int_{V}|\nabla \psi_{n}(x)|^{2} \leq \eta, \, \; \forall n \geq  N(\eta)
\end{equation}
provided that $a_{n}c_{n} \gg n^{3}$, i.e. $a_{n} n^{d-2} \gg n^{3}$
\end{prop}
{\it Proof.} Let us set:
\[ X_{n}:= \int_{V}|\nabla \psi_{n}(x)|^{2}  \]
For a given $\eta > 0$, we have :
\begin{align}
\label{wave35}
\E_{n} X_{n} & = \int_{0}^{\infty} \mu_{n}( X_{n} > y ) dy = \int_{0}^{2\eta} \mu_{n}( X_{n} > y ) dy  + \int_{2\eta}^{\infty} \mu_{n}( X_{n} > y ) dy\\ \nonumber
             & \leq 2\eta + \mu_{n}( X_{n} \leq \eta ) \int_{2\eta}^{\infty} \frac{\mu_{n}( X_{n} > y )}{\mu_{n}( X_{n} \leq \eta )} dy
\end{align}
Taking into account (\ref{wave25}) (or (\ref{wave26}) applied to $\int_{V}|\nabla \psi_{n}(x)|^{2} $), we will have:
\begin{align}
\label{wave36}
\int_{2\eta}^{\infty} \frac{\mu_{n}( X_{n} > y )}{\mu_{n}( X_{n} \leq \eta )} & \leq \int_{2\eta}^{\infty} e^{a_{n}n^{2} c_{n}(1-y /n^{2})} 
                                                                                   e^{- a_{n}n^{2} c_{n}(1-\eta /n^{2})} e^{\frac{+K |V| n^{3}}{\eta}} dy \\ \nonumber
																														& \leq \int_{2\eta}^{\infty} e^{-a_{n}n^{2} c_{n}(y-\eta) /n^{2}}e^{\frac{+K |V| n^{3}}{\eta}}dy \\ \nonumber	
																														& \leq     e^{\frac{+K |V| n^{3}}{\eta} - a_{n}n^{2} c_{n}\eta /n^{2}}
																														\int_{2\eta}^{\infty}  e^{-a_{n}n^{2} c_{n}(y-2\eta) /n^{2}} dy
\end{align}
The last integral is convergent (and is even converging to $0$) and the exponent 
$$ \exp [K |V| n^{3}/\eta - a_{n}n^{2} c_{n}\eta /n^{2}]$$
will converge to $0$ if $a_{n}n^{2} c_{n}\eta /n^{2} \geq 2 K |V| n^{3}/\eta $. For each fixed $\eta > 0$, this condition means that $a_{n}c_{n}\eta^{2} \geq 2 K |V| n^{3}$ and is satisfied when $a_{n}c_{n} \gg n^{3}$. So that under this condition the proposition follows from (\ref{wave35}) and (\ref{wave36}). $\Box$
\\
\\
\noindent
Proposition \ref{PropWave5} and lemma \ref{lemmaWave1} lead to the following:
\begin{co}
\label{CoroWave1}
For any $\varepsilon > 0 $, there exists a subsequence still denoted by $n$ and a volume $\tilde{V}_{\varepsilon}$ such that $|\tilde{V}_{\varepsilon}|\geq |V|(1-\sqrt{\varepsilon})$ and:
\begin{equation}
\label{wave40}
\E_{n}|\nabla \psi_{n}(x)|^{2} \leq \varepsilon, \, \; \forall n \geq  N(\varepsilon), \forall x \in \tilde{V}_{\varepsilon}
\end{equation}
provided that $a_{n}c_{n} \gg n^{3}$, i.e. $a_{n} n^{d-2} \gg n^{3}$. Furthermore, if we suppose that the field sequence $(\varphi, \mu_{n})$ has a weak limit then we also have:
\begin{equation}
\label{wave41}
\E_{n}|\psi_{n}(x)|^{2} \leq K \varepsilon, \; {\rm and} \; \E_{n}|\varphi_{n}(x)|^{2} \leq K \varepsilon c_{n} \; \forall n \geq  N(\varepsilon), \forall x \in \tilde{V}_{\varepsilon}
\end{equation}
for some constant $K$.
\end{co}
{\it Proof.} The first assertion is a consequence of lemma \ref{lemmaWave1} and proposition \ref{PropWave5}. Let us now show that (\ref{wave40}) implies also (\ref{wave41}). For notation convenience, let us take $x=0$, supposing that $0 \in \tilde{V}_{\varepsilon^{'}}$. We shall prove that $E_{n}|\psi_{n}(0)|^{2} \leq K \varepsilon$ which is equivalent to  $E_{n}|\varphi_{n}(0)|^{2} \leq K \varepsilon c_{n}$, a fact that will be valid for any other $x \in \tilde{V}_{\varepsilon^{'}} $. Let $f\in {\cal S}$ be a test function with support in $\tilde{V}_{\varepsilon^{'}}$. Then for any $x\in V$ we have :
\begin{equation}
\label{wave42}
\varphi_{n}(x)= \varphi_{n}(0)+ \int_{0}^{1}\nabla\varphi_{n}(tx).x dt
\end{equation}
and by theorem \ref{Distributions} (as we suppose that the field sequence $(\varphi, \mu_{n})$ has a weak limit):
\begin{equation}
\label{wave43}
\E_{\infty} \varphi(f)^{p} \sim \E_{n} \varphi(f)^{p} \sim \E_{n} \varphi_{n}(f)^{p}= \E_{n} (\int_{V} \varphi_{n}(x) f(x) dx)^{p}
\end{equation}
Hence:
\begin{align}
\label{wave44}
\E_{n} \varphi(f)^{2} \sim \E_{n} \varphi_{n}(f)^{2} & = \E_{n} \varphi_{n}(0)^{2} (\int_{V} f(x) dx)^{2}
                                 + \E_{n} [\int_{V} dx f(x) (\int_{0}^{1}\nabla\varphi_{n}(tx).x dt)]^{2} \\ \nonumber
                           &+ 2 \E_{n} \varphi_{n}(0) [\int_{V} dx f(x) (\int_{0}^{1}\nabla\varphi_{n}(tx).x dt)] \\ \nonumber
                    & \geq  (\int_{V} f)^{2} \E_{n} \varphi_{n}(0)^{2}  - |2 \E_{n} \varphi_{n}(0) [\int_{V} dx f(x) (\int_{0}^{1}\nabla\varphi_{n}(tx).x dt)]| \\ \nonumber
										& \geq (\int_{V} f)^{2} \E_{n} \varphi_{n}(0)^{2} - 2 (\E_{n} \varphi_{n}(0)^{2})^{1/2} \\ \nonumber
										& \times (\E_{n}  [\int_{V} dx f(x) (\int_{0}^{1}\nabla\varphi_{n}(tx).x dt)]^{2} )^{1/2}
\end{align}
Using again twice the Schwarz inequality we have : 
\begin{align*}
 [\int_{V} dx f(x) (\int_{0}^{1}\nabla\varphi_{n}(tx).x dt)]^{2} & \leq \int_{V} dx f(x)^{2} \int_{V} dx \int_{0}^{1}|\nabla\varphi_{n}(tx).x |^{2} dt 1_{f\neq 0} \\
                                                              & \leq \int_{V} dx f(x)^{2} \int_{V} dx \varepsilon^{2} c_{n}|x|^{2} dt 1_{f\neq 0} 
\end{align*}
and
\begin{equation}
\label{wave45}
\E_{n} \varphi(f)^{2} \geq (\int_{V} f)^{2} \E_{n} \varphi_{n}(0)^{2} - 2  \varepsilon^{2} c_{n} \int_{V} dx f(x)^{2} \int_{V}|x|^{2} dt 1_{f\neq 0} dx (\E_{n} \varphi_{n}(0)^{2})^{1/2} 
\end{equation}
Since $\E_{n} \varphi(f)^{2} < \infty$ and $\varphi_{n}(0)^{2}= c_{n}\psi_{n}(0)^{2}$, (\ref{wave41}) follows from (\ref{wave45}) for $x=0$. This is also valid for all $x \in \tilde{V}$. We also note that the constant $K$ does not depend on the volume $|V|$. $\Box$
\\
\noindent
{\it Proof of Proposition \ref{PropWave1b}.} By Corollary \ref{CoroWave1} we have the estimate $\E_{n}\varphi_{n}(x)^{2} \leq K \varepsilon c_{n}$ for $n\geq n(\varepsilon)$ which is stronger than the inequality (\ref{wave9}). Taking into account this improvement, the proof of Proposition \ref{PropWave1b} as exactly the same than that of Proposition \ref{PropWave1}. $\Box$
\\
\noindent
{\it Proof of Proposition \ref{PropWave1c}.}
We shall outline the main steps of this proof. We start from the result we have just obtained in the most general case $a_{n} n^{2} \gg g_{n}c_{n}$, that is :
\begin{equation}
\label{wave46}
\E_{n} \varphi_{n}^{2}(x) \ll c_{n} \sim K n^{d-2}
\end{equation}
If we suppose that $\mu_{n} \rightarrow \mu_{\infty}$ weakly, then one may expect that we will also have:
\begin{equation}
\label{wave47}
\E_{\infty} \varphi_{n}^{2}(x) \ll c_{n} \sim K n^{d-2}
\end{equation}
If this is true, then this will imply that the $\varphi^{4}_{d}$ will be trivial and even identically vanishing: $\varphi\equiv 0$. This is a consequence of a general property of the quantum fields that we briefly recall: For the free field of covariance $C$, we know that:
\begin{equation}
\label{wave48}
C(x-y) \sim \frac{K}{|x-y|^{d-2}}
\end{equation}
and a quantum Euclidean field with covariance kernel $C_{\nu}$ is said to have a canonical short-distance behaviour if (\ref{wave48}) is still satisfied by $C_{\nu}$. Now it is easy to see that this condition implies also that $\E_{\nu} \varphi_{n}^{2}(x) \sim K c_{n} \sim K n^{d-2}$. On the other hand, we also have the following fact (see \cite{FFS} p.396): If an Euclidean quantum field has a covariance $C$ less singular than that of the free field and satisfies the Osterwalder-Schraeder axioms, then its 2-point function is identically null ($S^{(2)}\equiv 0$), in particular when this field is given by a measure $\nu$ in ${\cal S}^{'}(\R^{d})$, we will have $\E_{\nu} \varphi_{n}^{2}(f)=0 $ for all $f\in {\cal S}(\R^{d}) $ which means that $\nu=\delta_{0}$. This a consequence of The K\"allen-Lehmann representation theorem. 
In our case  if $\E_{\infty} \varphi_{n}^{2}(x) \ll c_{n} \sim K n^{d-2}$ then the field corresponding to $\mu_{\infty}$ will necessarily be less singular than the free field and therefore we will have $\E_{\mu_{\infty}} \varphi_{n}^{2}(f) $ for all $f\in {\cal S}(\R^{d}) $ and $\mu_{\infty}=\delta_{0}$.

It remains to see whether (\ref{wave47}) is valid. We first note that (\ref{wave47}) is not a trivial consequence of (\ref{wave46}). To prove it, we are led to consider the double sequence $\E_{n} \varphi_{k}^{2}(x)$ and in order to determine an estimate of it, we consider the following action:
\begin{align}
\label{wave49}
{\cal A}_{n,k}(\varphi) & = \int_{V}  dx [ g_{n}: \varphi^{4}_{k}:+ m_{n} : \varphi^{2}_{k}: + a_{n}: (\partial \varphi_{k})^{2}:] \\ \nonumber
                    & = a_{n}c_{n}n^{2} {\cal A}_{n}^{(1)} (\psi_{k}) 
\end{align}
with
\begin{equation}
\label{wave50}
{\cal A}_{n}^{(1)} (\psi_{k}) = \int_{V} [ \frac{ g_{n} c_{n}}{a_{n}n^{2}} :\psi_{k}^{4}:+ \frac{m_{n}}{a_{n}n^{2}} : \psi_{k}^{2}: 
                                                + \frac{1}{d_{n}} \int_{V}: (\partial\psi_{k})^{2}:]
\end{equation}
and the corresponding measures:
\begin{equation}
\label{wave51}
d\mu_{n,k}(\varphi) = \frac{e^{-{\cal A}_{n,k}(\varphi_{k})}}{Z_{n,k}} = \frac{e^{-a_{n}c_{n}n^{2}{\cal A}_{n,k}(\psi_{k})}}{Z_{n,k}}
\end{equation}
Now it can be shown that all the previous results (beginning with the large deviations lower bound) are valid for the sequence $\mu_{n,k}$ when we take a double limit $n, k \rightarrow \infty, k< n$. As a result we will get:
\begin{equation}
\label{wave52}
\lim_{n,k \rightarrow \infty, k< n} \E_{n,k}(\partial \psi)(x)^{2} =0
\end{equation}
and $\lim_{k \rightarrow \infty} \E_{\infty} (\partial \psi_{k})(x)^{2} =\lim_{k \rightarrow \infty} \E_{\infty} \psi_{k}(x)^{2}= 0$ which yields 
$\E_{\infty} \varphi_{k}(x)^{2} \ll c_{n}$. On the other hand, one can verify that the previous estimates are also independent of the volume and that (\ref{wave47}) is also valid as we take the infinite volume. This completes the outline of the proof of Proposition \ref{PropWave1c}.$\Box$

\begin{Rk}
It may be possible to have a more direct proof of Proposition \ref{PropWave1c}. However such a direct proof under the most general condition $a_{n} n^{2} \gg g_{n}c_{n}$ seems not to be easy. While Propositions \ref{PropWave1} and especially \ref{PropWave1b} provides the triviality result $\Delta \varphi=0 $ under fairly general conditions,
the proof that $\E_{\infty} \varphi(f)^{2}= 0$ does not seem to be straightforward. Nevertheless, this can be shown under a little more restrictive conditions: at first, we use the dynamic equation $\E_{n} \varphi(f) \varphi(\tilde{h})= (h,f) - \E_{n} \varphi(f) \partial_{h} \A_{n}$ which we write in the form :
\begin{align}
\label{wave53}
\E_{n} \varphi(f) \varphi(\tilde{h}) &= (h,f) -  4g_{n} \E_{n} \varphi (f) \int_{V}\varphi_{n}^{3}(x) h_{n}(x) + 4g_{n}\times 3c_{n}  \E_{n} \varphi (f) \varphi_{n}(h_{n}) \\ \nonumber
                                & - 2m_{n} \E_{n} \varphi (f) \varphi_{n}(h_{n})  + 2a_{n} \varphi (f) \E_{n}\varphi_{n}(\Delta h_{n})
\end{align}
The idea is that since the terms $ \E_{n}\varphi_{n}^{2}(x)$ is bounded by $K\varepsilon c_{n}$ by Corollary \ref{CoroWave1}, the dominant term in (\ref{wave53}) will be $12 g_{n}c_{n}  \E_{n} \varphi (f) \varphi_{n}(h_{n}) $ if $a_{n}$ is not too large and (\ref{wave53}) will give $\E_{n} \varphi (f) \varphi_{n}(h_{n}) \rightarrow 0$. More precisely, we begin first by controlling the term containing $\varphi_{n}^{3}(x)$: By Griffiths and Gaussian inequalities (recall that $f\geq 0$) we have:
\begin{equation}
\label{wave54}
   0 \leq \E_{n} \varphi (f)\varphi_{n}^{3}(x) \leq 3 \E_{n} \varphi (f) \varphi_{n}(x) \E_{n} \varphi_{n}^{2}(x) \leq 3 \varepsilon c_{n}  \E_{n} \varphi (f) \varphi_{n}(x)
\end{equation}
Therefore (\ref{wave53}) yields:
\begin{align}
\label{wave55}
12 g_{n}c_{n}  \E_{n} \varphi (f) \varphi_{n}(h_{n})  & \leq \E_{n} \varphi(f) \varphi(\tilde{h})-(h,f) + 4g \times 3 \varepsilon c_{n}  \E_{n} \varphi (f) \varphi_{n}(x) \\ \nonumber
                                      & + 2m_{n} \E_{n} \varphi (f) \varphi_{n}(h_{n}) 
																			+ 2a_{n} |\E_{n}\varphi_{n}(\Delta h_{n})|
\end{align}
This inequality shows that if $g_{n}c_{n} \geq K a_{n}$ and $ g_{n}c_{n} \gg m_{n} $ we will have 
\begin{equation}
\label{wave56}
\E_{n} \varphi (f) \varphi_{n}(h_{n}) \rightarrow 0.
\end{equation}
As we work in the case $  a_{n} n^{2} \gg g_{n}c_{n} $ with $ c_{n} \sim K n^{d-2}$, we deduce that (\ref{wave49}) holds when :
\begin{equation}
\label{wave57}
            g_{n} n^{d-4}  \ll a_{n}  \leq K'  g_{n} n^{d-2}  \; \; {\rm and } \; \; g_{n} n^{d-2} \gg m_{n}
\end{equation}
On the other hand (\ref{wave53}) yields also:
\begin{align}
\label{wave58}
2m_{n} \E_{n} \varphi (f) \varphi_{n}(h_{n}) & \leq \E_{n} \varphi(f) \varphi(\tilde{h})-(h,f) + 4g \times 3 \varepsilon c_{n}  \E_{n} \varphi (f) \varphi_{n}(x) \\ \nonumber
                                      & + 12 g_{n}c_{n}  \E_{n} \varphi (f) \varphi_{n}(h_{n})
																			+ 2a_{n} |\E_{n}\varphi_{n}(\Delta h_{n})|
\end{align}
and (\ref{wave56}) is valid when $m_{n} \gg g_{n}c_{n}$ and $   m_{n} \gg a_{n} $ which was already stated in $\S 5.1.$ but is recalled here as a particular case when the mass term is predominant ($m_{n}\ll a_{n} n^{2}$).
\end{Rk}
\subsection{The case where the interaction term is predominant : $g_{n}c_{n} \gg m_{n}, a_{n} n^{2}$}

In this case we write the action as:
\begin{align}
\label{int3}
{\cal A}_{n}(\varphi)= & g_{n} \int_{V}: \varphi_{n}^{4}:+ m_{n} \int_{V}: \varphi_{n}^{2}:+a_{n} \int_{V}: (\partial\varphi_{n})^{2}:  \nonumber \\
            = & g_{n} c_{n}^{2} \int_{V}: \psi_{n}^{4}:+ m_{n}c_{n} \int_{V}: \psi_{n}^{2}: 
						       + a_{n} c_{n} d_{n} \frac{1}{d_{n}} \int_{V}: (\partial\psi_{n})^{2}:  \\ \nonumber
						= & g_{n} c_{n}^{2} {\cal A}_{n}^{(1)} (\varphi) \\ \nonumber
\end{align}
with:
\begin{equation}
\label{Int4}
{\cal A}_{n}^{(1)} (\varphi)= \lambda_{n}\int_{V}: \psi_{n}^{4}:+ \alpha_{n} \int_{V}: \psi_{n}^{2}: 
						       +  \frac{\beta_{n}}{d_{n}} \int_{V}: (\partial\psi_{n})^{2}: 
\end{equation}
\begin{equation}
\label{Int5}
\psi_{n} := \varphi/\sqrt{c_{n}}, \; \lambda_{n} = 1, \; \alpha_{n}=\frac{m_{n}}{g_{n}c_{n}}, \; \beta_{n}=\frac{a_{n}d_{n}}{g_{n}c_{n}}, \; d_{n}=n^{2}.
\end{equation}
As we can work with a subsequence, we suppose, without restriction, that $\alpha_{n}\rightarrow \alpha$ and  $\beta_{n}\rightarrow \beta$, with $\alpha, \beta < +\infty$. The principle of the least action leads to consider the classical action :
\begin{equation}
\label{Int6}
{\cal A}_{n}^{c} (\chi)= \int_{V} H(\chi(x)):= \int_{V}: \chi^{4}:+ \alpha \int_{V}: \chi^{2}: 
						       + \beta \frac{1}{d_{n}} \int_{V}: (\partial\chi)^{2}: 
\end{equation}
It is understood here that $\chi : \R^{d}\rightarrow \R$ are deterministic functions and we note symbolically $: \chi^{4}:= \chi^{4} - 6 \chi^{2} + 3$ and $:\chi^{2}:=\chi^{2}-1$.

The minimum of the action ${\cal A}_{n}^{(1)} (\chi)$ is attained for $\chi \equiv constant =\chi_{0}$ with $4 \chi_{0}^{3} -(6-\alpha) \chi_{0}=0 $, which leads to the following cases:
\\
\noindent
$\bullet$ Case 1.1 : $\alpha < 6 $ and $\chi_{0}^{2}= 3 - \alpha/2$, i.e. the minimum of the action is attained for $\chi_{0}= \pm \sqrt{ 3 - \alpha/2}$
\\
\noindent
$\bullet$ Case 1.2 : $\alpha = 6 $ in which case $H(X) = X^{4}-3$ and the minimum is attained for $X=0$ i.e. $\chi_{0}= 0$.
\\
\noindent
$\bullet$ Case 1.3 : $\alpha > 6 $ which means that $m_{n}/g_{n}c_{n} > 6$ for large $n$, this case was considered in section 5.1.
\\
\noindent
For the case 1.1 we have the following result:
\begin{prop}
\label{PropIF2}
In the case 1.1, in particular when : $ g_{n} c_{n} \gg m_{n} $ or $m_{n}\sim \alpha g_{n}c_{n}, \alpha < 6 $  and $ g_{n} c_{n} \gg  n^2 a_{n}$, the regularized fields $\varphi_{n}$ and $\psi_{n}=\varphi_{n}/\sqrt{c_{n}}$ may be written as :
\begin{equation}
\label{Int7}
\psi_{n}(x) = \sqrt{3} s_{n}(x) + \xi_{n}(x), \; \; \varphi_{n}(x) = \sqrt{c_{n}}(\sqrt{3-\alpha/2} s_{n}(x) + \xi_{n}(x))
\end{equation}
where $s_{n}(x), \xi_{n}(x)$ are random fields satisfying:
\begin{equation}
\label{Int7}
s_{n}(x) = \mp 1,
\end{equation}
and for any $\varepsilon >0 $, there exists $n_{\varepsilon}\geq 1$ such that for all $n \geq l$: 
\begin{equation}
\label{Int7}
\int_{V} \xi_{n}(x)^{2} dx <  \varepsilon \; \; {\rm and} \; \;   \E_{n}\xi_{n}(x)^{2}  <  \varepsilon \; \forall x \in V_{l}
\end{equation}
where the sequence of volumes $V_{l}$ verifies: $V_{l} \uparrow V$ and $|V_{l}| \geq (1-2\sqrt{\varepsilon}) |V|$.
\end{prop}
{\it Proof.} 
Whe prove the proposition when $\alpha=0$, i.e. $g_{n}c_{n} \gg m_{n}, a_{n}n^{2} $. The proof is similar if $ 0<\alpha < 6 $.  The minimum of the functional ${\cal A}_{n}^{c}$ (\ref{Int6}) is $6$, and the principle of the least action tells us that for any $\varepsilon  > 0 $, there is some $n_{\varepsilon}\geq 1$ such that:
\begin{equation}
\label{Int8}
-{\cal A}_{n}^{(1)} (\varphi)= - [\lambda_{n}\int_{V}: \psi_{n}^{4}:+ \alpha_{n} \int_{V}: \psi_{n}^{2}: 
						       + \beta_{n} \frac{1}{d_{n}} \int_{V}: (\partial\psi_{n})^{2}: ]\geq 6 -\varepsilon 
\end{equation}
Since $\alpha_{n}\sim 0, \beta_{n}\sim 0$, this implies that 
\begin{equation}
\label{Int}
-\int_{V}: \psi_{n}^{4}: \geq 6 -\varepsilon , \; \forall n\geq n_{2, \varepsilon}
\end{equation}
 for some $n_{2, \varepsilon}$. If we set:
\[ \psi_{n}^{2}(x)= 3+ \varepsilon_{n}(x) \]
then (\ref{Int9}) implies that for $n\geq n_{2, \varepsilon}$:
\[ \int_{V}\varepsilon_{n}^{2}(x) \leq \varepsilon \]
This implies of course that $\int_{V} \E_{n}\varepsilon_{n}^{2}(x) \leq \varepsilon$. By lemma \ref{lemmaWave1}, there is a subsequence that we also denote by $n$, such that for any $\varepsilon>0$, there exists $l\geq 1 $: $\E_{n}\varepsilon_{n}^{2}(x) \leq \varepsilon$ for all $x\in V_{l}$ and $n\geq l$ with a sequence of increasing volumes $V_{l}$ that satisfy the conditions of the proposition. $\Box$
\\
\noindent
Proposition \ref{PropIF2} means that the normalized field $\psi_{n}$ is rapidly oscillating between the values $-\sqrt{3}$ and $+\sqrt{3}$ when $n\rightarrow \infty$. If $\alpha \neq 0$, i.e. $m_{n}\sim \alpha g_{n}c_{n}, \alpha < 6 $, $\psi_{n}$ will oscillate between the values $-\sqrt{3-\alpha/2}$ and $+\sqrt{3-\alpha/2}$.

\end{document}